\documentclass[11pt,reqno]{amsart}
\usepackage{a4wide}
\usepackage{float,caption,capt-of,scalefnt,relsize}
\usepackage{euscript,amsmath,amssymb,amsbsy}
\usepackage{array,longtable,tabularx,multirow,enumitem,url}
\usepackage{graphicx,float}\usepackage{subfigure,epsfig}

\usepackage[colorlinks=true,linkcolor=black,urlcolor=black,citecolor=black]{hyperref}

\numberwithin{equation}{section}\numberwithin{figure}{section}

\newcounter{msct}[section]\renewcommand{\themsct}{\thesection.\arabic{msct}}

\newenvironment{m-theorem}{\vskip5pt\refstepcounter{msct}\trivlist \itemindent 0pt%
\item[\hskip\labelsep\bf Theorem~\themsct]\it\ignorespaces}{\endtrivlist\vskip3pt}

\newenvironment{m-proposition}{\vskip5pt\refstepcounter{msct}\trivlist \itemindent0pt%
\item[\hskip\labelsep\bf Proposition~\themsct]\it\ignorespaces}{\endtrivlist\vskip3pt}

\newenvironment{m-corollary}{\vskip5pt\refstepcounter{msct}\trivlist \itemindent 0pt%
\item[\hskip\labelsep\bf Corollary~\themsct]\it\ignorespaces}{\endtrivlist\vskip3pt}

\newenvironment{m-lemma}{\vskip5pt\refstepcounter{msct}\trivlist \itemindent 0pt%
\item[\hskip\labelsep\bf Lemma~\themsct]\it\ignorespaces}{\endtrivlist\vskip3pt}

\newenvironment{m-definition}{\vskip5pt\refstepcounter{msct}\trivlist \itemindent0pt%
\item[\hskip\labelsep\bf Definition~\themsct]\ignorespaces}{\endtrivlist\vskip5pt}

\newenvironment{m-notation}{\vskip5pt\refstepcounter{msct}\trivlist \itemindent0pt%
\item[\hskip\labelsep\bf Notation~\themsct]\ignorespaces}{\endtrivlist\vskip5pt}

\newenvironment{m-example}{\vskip5pt\refstepcounter{msct}\trivlist \itemindent0pt%
\item[\hskip\labelsep\bf Example~\themsct]\ignorespaces}{\endtrivlist\vskip5pt}

\newenvironment{m-remark}{\vskip5pt\refstepcounter{msct}\trivlist \itemindent0pt%
\item[\hskip\labelsep\bf Remark~\themsct]\ignorespaces}{\endtrivlist\vskip5pt}

\newenvironment{m-procedure}{\vskip5pt\refstepcounter{msct}\trivlist\itemindent0pt%
\item[\hskip\labelsep\bf Procedure~\themsct]\ignorespaces}{\hfill\endtrivlist\vskip5pt}%

\newenvironment{m-question}{\vsk\char{cmmib10}{15}ip5pt\refstepcounter{msct}\trivlist \itemindent0pt%
\item[\hskip\labelsep\bf Question.]\ignorespaces}{\endtrivlist\vskip5pt}

\newenvironment{thm-nono}[1]{\vskip5pt\trivlist \itemindent 0pt %
\item[\hskip\labelsep\bf Theorem~{\rm\mbox{#1}}]\it\ignorespaces}{\endtrivlist\vskip5pt}

\newenvironment{lm-nono}[1]{\vskip5pt\trivlist \itemindent0pt%
\item[\hskip\labelsep\bf Lemma~{\rm\mbox{#1}}]\it\ignorespaces}{\endtrivlist\vskip5pt}

\newenvironment{conj-nono}[1]{\vskip5pt\trivlist \itemindent0pt%
\item[\hskip\labelsep\bf Conjecture~{\rm\mbox{#1}}]\it\ignorespaces}{\endtrivlist\vskip5pt}

\newenvironment{m-thank}{\vskip5pt\trivlist \itemindent0pt%
\item[\hskip\labelsep\it Acknowledgments]\ignorespaces}{\endtrivlist\vskip5pt}

\newenvironment{m-proof}{\vskip2pt\trivlist \itemindent0pt%
\item[\hskip\labelsep\it Proof.]\ignorespaces}{\hfill$\Box$\endtrivlist\vskip5pt}%

\newenvironment{m-asmp}{\vskip5pt\trivlist \itemindent0pt%
\item[\hskip\labelsep\bf Assumption.]\ignorespaces}{\hfill\endtrivlist\vskip5pt}%

\newcounter{meqn}[section]\renewcommand{\themeqn}{\thesection.\arabic{meqn}}
\newenvironment{m-eqn}[1]{\vskip5pt\refstepcounter{meqn}%
\trivlist\itemindent0pt\item[]\ignorespaces%
\hfill $\displaystyle #1$\hfill\hbox{\rm(\themeqn)}}{\endtrivlist\vskip5pt}

\newcommand{\bibauth}[2]{\textrm{{#1}~{#2},}}
\newcommand{\bibtitl}[1]{\textit{#1}.}
\newcommand{\bibjnyp}[4]{\textrm{#1} \textbf{#2} (#3), {#4}.}
\newcommand{\bibinbook}[4]{In: \textrm{#1}\textrm{, #2}\textrm{, #3}\textrm{, #4}.}

\let\mt\mapsto

\let\mbb\mathbb

\DeclareFontFamily{OT1}{rsfs}{}
\DeclareFontShape{OT1}{rsfs}{n}{it}{<->rsfs10}{}
\DeclareMathAlphabet{\crl}{OT1}{rsfs}{n}{it}

\let\disp\displaystyle
\let\ges\geqslant

\let\dta\delta 
\let\les\leqslant
\let\mt\mapsto
\let\nit\noindent

\newcommand\uset[2]{\mathop{#2}\limits_{#1}}

\newcommand{\ala}{\alpha}
\newcommand{\aaa}{a}
\let\bta\beta
\newcommand{\bvp}{BVP}\newcommand{\ivp}{IVP}
\newcommand{\bbb}{b}

\newcommand{\cc}{c}\newcommand{\CC}{C}
\newcommand{\ct}{{\it ct}}\newcommand{\CT}{{\it CT}}
\newcommand{\ee}{{e}}\newcommand{\re}{{\rm e}}
\newcommand{\exc}{{\rm ex}}
\newcommand{\kk}{{k}}\newcommand{\KK}{{K}}
\let\kpa\kappa
\let\Lda\Lambda
\newcommand{\LL}{L}
\newcommand{\ave}{{av}}
\newcommand{\nn}{{n}}\newcommand{\NN}{{N}}
\newcommand{\ode}{ODE}
\newcommand{\pp}{{p}}
\newcommand{\qq}{{q}}
\newcommand{\resd}{{\rm res}}

\newcommand{\uu}{u}
\newcommand{\vv}{v}
\newcommand{\ww}{w}\newcommand{\WW}{W}
\newcommand{\yy}{y}
\newcommand{\zz}{z}

\let\veps\varepsilon


\begin{document}

\title[Troesch's equation with charge separation]{Troesch's equation with charge separation}
\author{Mihai Halic} 

\begin{abstract}
We consider the following generalization of Troesch's classical two-point-BVP in order to take into account charge separation, which is required by physical considerations: $w''=L\sinh(Lw), w(0)=W_0<0<W_1=w(1)$. In contrast with the classical case, the solution of this generalized equation possesses boundary layer at both ends. We determine an explicit, approximate solution and verify its precision in several ways. 

The analytic solution hinges on a result concerning a class of overdetermined first order BVPs. Presenting this matter constitutes the second objective of the article. 
\end{abstract}

\maketitle

\section*{Introduction}

This work is rooted in Weibel's articles~\cite{weib1,weib2}, devoted to studying the motion of plasma in magnetic field. Troesch~\cite{tsch} simplified Weibel's system of differential equations, leading him to the equation, 
$$
\yy''={\LL}\sinh(\LL\yy),\;\yy(0)=0,\;\yy(1)=1\;\;(\LL>0), 
$$
which is the subject of considerable scholarly literature. Our starting point, that will be explained, is that a simplifying assumption leading to this boundary value problem ({\bvp}, for short) neglects an essential aspect in Weibel's work: \textit{charge separation}. 
For this reason, we return to his original work and introduce another second order {\bvp} ---we call it `generalized Troesch-problem'--- which takes into account this matter:
\begin{align*}
\yy''=\frac{\LL}{2}\cdot\bigl(\qq^{-1}\re^{\LL\yy}-\qq\re^{-\LL\yy}\bigr),\;\;\yy(0)=0,\;\yy(1)=1\;(\LL>0,\;\qq>1).\tag{$\ast$}\label{eq:star}
\end{align*}
This modification is not only a mathematical artifice: for sufficiently large values $\LL$, the solution possesses \emph{two boundary layers}, at the left- and right-end, in contrast with the classical Troesch equation. The double boundary layer makes trickier to understand the behaviour of the solution: saying that `the initial derivative is almost zero' is no longer true. Instead, the solution is almost zero about some unknown point between $0$ and $1$.

The equation can't be explicitly solved and the parameters make more challenging to understand the dependence of the solution on them. This motivates the necessity to develop an analytical method for approximately solving this  (and similar) equation(s). 

The article is structured as follows. In Section~\ref{sct:ode}, we explain how autonomous, second order {\bvp}s lead to an (artificially looking) type of transcendental equations, whose (approximate) analytic solution is derived in the Appendix. We chose to postpone the matter because it's slightly tedious and we wished to prioritize the application; however, this `toolbox' is indispensable. The application ---that is, the generalized Troesch-problem--- is considered in the second section. 
\begin{itemize}[leftmargin=3ex]
\item 
We show how charge separation leads to \eqref{eq:star}. For $\LL$ large compared to $\qq$, its solution possesses two boundary layers (near $0$ and $1$). 
\item 
An unexpected `superposition principle' allows breaking the strongly non-linear equation~\eqref{eq:star} into two independent (decoupled) Troesch-type equations (with modified right-boundary condition).
\item 
We use the `analytical toolbox' to produce explicit, approximate solutions of the decoupled {\bvp}s. Their difference is our approximate solution of \eqref{eq:star}.
\item 
We verify the precision of this approximate solution for several parameters. It turns out that it's very precise, the error decreases exponentially with $\LL$. 
\end{itemize}


\section{Approximate solutions of autonomous, second order {\bvp}s}\label{sct:ode}

Yamakawa-Kreinovich~\cite{ya-kr} explain the important role played by second-order {\ode}s in physics. Fundamental phenomena are modelled either by first or second order differential equations~[\textit{id}., p.\,1765]: \emph{``Every time we have a physical process whose evolution is described by a function of three or more variables, this process is not fundamental.''} 
This justifies the interest for analytical methods to determine approximate solutions of two-point {\bvp}s: 
\begin{align*}
\yy''=f(\yy),\;\;\yy(\aaa)=\yy_\aaa,\;\yy(\bbb)=\yy_\bbb.
\tag{BVP2}\label{eq:bvp2}
\end{align*} 
Numerical methods are convenient and accurate as long as the solution $\yy$ doesn't have large variation. In `unpleasant cases', the solution has a very thin boundary layer. That is, the graph of $\yy$ is almost vertical and implementing algorithms becomes difficult, often leading to indeterminacies. 
Conceptually, numerical methods give no rigorous insight into the dependence on parameters ($\yy_a, \yy_b$, or $f$ itself); each data requires its own computations. 

A different approach is to multiply the equation by $\yy'$ and integrate it:
\begin{align*}
\yy'=\pm\sqrt{2}\cdot\sqrt{F(\yy)+\dta},\;\;\yy(\aaa)=\yy_\aaa,\;\yy(\bbb)=\yy_\bbb,\;\dta\in\mbb R,
\tag{BVP1}\label{eq:bvp1}
\end{align*}
where $F$ is an anti-derivative of  $f$. This \textit{first order} problem is overdetermined: the value of unknown parameter $\dta$ is determined by the boundary conditions, the problem is solvable only for a specific value. The difficulty is that, except in trivial situations, the value of $\dta$ is impossible to compute. Even approximating it is delicate: to initiate the shooting-method, one needs an initial input and this may be very challenging. 

All together, these facts prompted us to develop analytical approach to~\eqref{eq:bvp1}. It involves two stages: changing it into an integrable first-order {\ode} (straightforward, below) and approximating $\dta$ (somewhat tedious, Appendix~\ref{sct:back}). The reason for deferring the latter issue to the Appendix is to reach faster the application to Weibel and Troesch's equation in \S\ref{sct:aplic}. 


\subsection*{Algorithm for approximating solutions of (BVP1)}

We replace the initial differential equation with an integrable one. Its solution approximates the solution of the given problem. 

\begin{enumerate}[leftmargin=5ex]
\item 
Computational difficulties arise when the equation admits boundary layer. Often it's due to a dominant term $Z(y)$ contained in $f(y), F(y)$. The change of variables $z=Z(y)$ simplifies the problem, the equation in $z$ has polynomial growth:
\begin{m-eqn}{
z'=K\cdot\sqrt{F(z)+\dta},\;\;F\;\text{has polynomial groth}.
}\label{eq:zF}
\end{m-eqn}
\item 
We replace the right-hand-side above by a (polynomial) truncated Taylor series $P(z,\dta)$. Note that 
$\,z'=P(z,\dta)\,$ is \emph{explicitly integrable only} for $P$ of degree \emph{at most two} in $\zz$. Thus we must take quadratic expansions of $\sqrt{F(z)+\dta}$. 
\item 
We reach the equation~\eqref{eq:bvp} in the Appendix, estimate $\dta$ by using the algorithm there, and insert it into $P(z,\dta)$.  
The new differential equation $\zz'=P(\zz,\dta)$ \emph{approximates} the initial one and is integrable, we obtain an approximate solution of~\eqref{eq:zF}. (The accuracy of the approximation hinges on the precision of the estimate for $\dta$. This explains the tedious work in the Appendix.) The approximate solution of~\eqref{eq:bvp2} is obtained by substituting $y=Y(z)$, $Y=Z^{-1}$. 
\end{enumerate}


\section{Generalizing the Troesch-equation}\label{sct:aplic}

The motivation stems from Weibel's work~\cite{weib1,weib2}, investigating the dynamics of plasma in a magnetic field. The process is described by a system of differential equations, difficult to solve numerically. A combination of these laws (\cite[eq. (16, 21, 27, 28)]{weib1}, \cite[Appendix]{tsch}) yields:
$$
\uu''=\frac{\kpa}{2}\cdot(\NN_i-\nn_\ee),\;\;
\NN_i=\NN_0\cdot\re^{\KK\uu-\CT\cdot E_0^2},\;\nn_\ee=\nn_0\cdot \re^{-\KK\uu-\ct\cdot E_0^2}.
$$ 
(The symbols $\kpa,\KK,\CT,\ct$ stand for (positive) physical quantities.) The field-strength $E_0$ obeys a second order, non-linear ODE (which is satisfied by the zero-function). The number of ions (resp. electrons) in the plasma are given by the functions $\NN_i, \nn_\ee$. 

Below we reproduced the expected shapes of the graphs of $\NN_i, \nn_\ee$. Weibel's work was revisited by Braga, Soares, \textit{et al.}~\cite{br-so, so-br} using modern computational power. They found that Weibel's plots occur for a low number of grid points; finer grids lead to modified electronic density shapes. 

In any event, there is \emph{`charge separation'}: $(\nn_\ee-\NN_i)$ has opposite signs at the endpoints of the interval of definition. This motivates our considerations: our aim is to modify Troesch's {\bvp}, so that it takes into account the necessary charge separation. 
\\[1ex]\begin{minipage}[b]{.95\textwidth}\centering
\begin{minipage}[b]{.45\textwidth}
\begin{figure}[H]
\includegraphics[height=.7\textwidth]{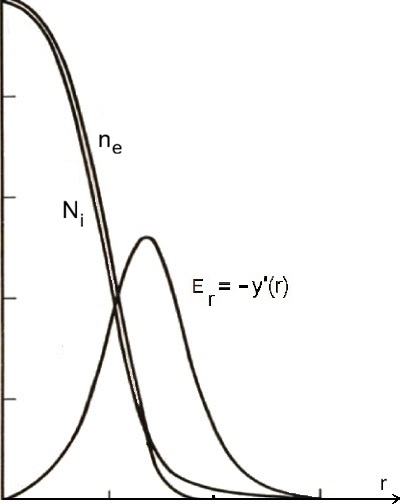}
\end{figure}
\end{minipage}
\begin{minipage}[b]{.45\textwidth}
\includegraphics[height=.73\textwidth]{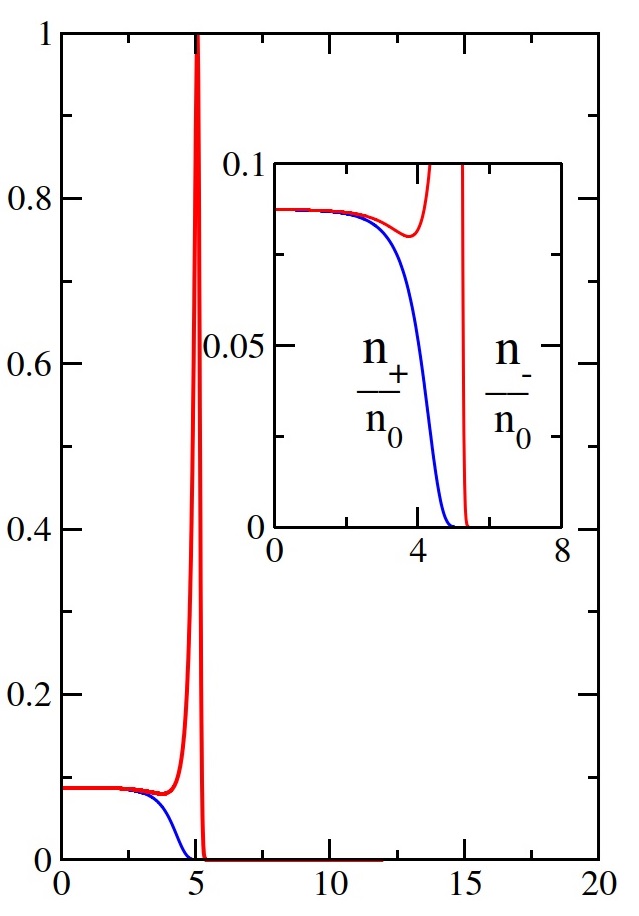}
\end{minipage}
\captionof{figure}{graphs of $N_i$ and $n_\ee$ (cf.~\cite[Fig. 3]{weib1} left, and~\cite[Fig. 3]{br-so} right).}
\label{fig:weibel}
\end{minipage}


\subsection{Setting up the problem}\label{ssct:bvp}

For mathematical convenience, to deal with increasing functions, we change the variable $r=1-x$. 
The first simplification of the original equations is to assume $E_0=0$ or small, so it can be neglected or absorbed into $\kpa$. We get 
$$\uu''=\frac{\kpa}{2}\cdot(N_0\re^{\KK\uu}-n_0\re^{-\KK\uu})
\text{ and }\NN=\NN_0\re^{\KK\uu}, \nn=\nn_0\re^{-\KK\uu}.$$ 
The substitutions $\uu:=\sqrt{\frac{\kpa\sqrt{\NN_0\nn_0}}{\KK}}\yy, \LL:=\sqrt{\kpa\KK\sqrt{\NN_0\nn_0}},\qq:=\sqrt{\frac{\nn_0}{\NN_0}}$ yield the `generalized' Troesch equation
\begin{m-eqn}{
\begin{array}{lll}
\disp\yy''=\frac{\LL}{2}\big(\qq^{-1}\re^{\LL\yy}-\qq\re^{-\LL\yy}\big),
\\[1ex] 
\yy(0)=0, \yy(1)=1,&&\text{(imposed boundary conditions)}.
\end{array}
}\label{eq:gentrs-v1}
\end{m-eqn}
Charge separation translates into: 
\begin{m-eqn}{
\{\yy''(0)<0,\; \yy''(1)>0\}\;\Leftrightarrow\;1<\qq<\re^{\LL}.
}\label{eq:qL}
\end{m-eqn}
The equation can be brought into a more familiar form. For $\pp:=\frac{\ln(\qq)}{\LL}, \ww:=\yy-\pp,$ we obtain:
\begin{m-eqn}{
\ww''=\LL\cdot\sinh(\LL\ww),\;\, \ww(0)=\WW_0,\;\ww(1)=\WW_1,\;\,(\WW_0<0<\WW_1).
}\label{eq:gentrs-v2}
\end{m-eqn}
We denote by $\ww_\exc$ the exact solution and consider arbitrary $\WW_0<0<\WW_1$. \emph{To our knowledge, \eqref{eq:gentrs-v2} has not been investigated so far.} Note that physically relevant values $\LL$ are usually large. Also, \eqref{eq:qL} is equivalent to $\pp\in(0,1)$; for `physically motivated' cases, one should think off  
\begin{m-eqn}{
\WW_1=1-\pp\in(0,1),\;\WW_0=-\pp\in(-1,0).
}\label{eq:Wp}
\end{m-eqn}
For $\NN_0=\nn_0$ one recovers the `classical' Troesch's {\bvp}; however, it neglects charge separation. The defining {\ode} of~\eqref{eq:gentrs-v2} is the same, but has modified boundary values. This change determines a very different behaviour of the solution. 

The author determined~\cite{csfx} very precise approximate analytical solutions of the classical Troesch equation, going far beyond the reach of pure computational approaches. In the current situation, applying \textit{ad litteram} the same approach fails badly: the technique doesn't handle two boundary layers. We start with a few plots. 
\\[1ex]\begin{minipage}{.99\textwidth}\centering 
\begin{minipage}[b]{.25\textwidth}
\includegraphics[height=.9\textwidth]{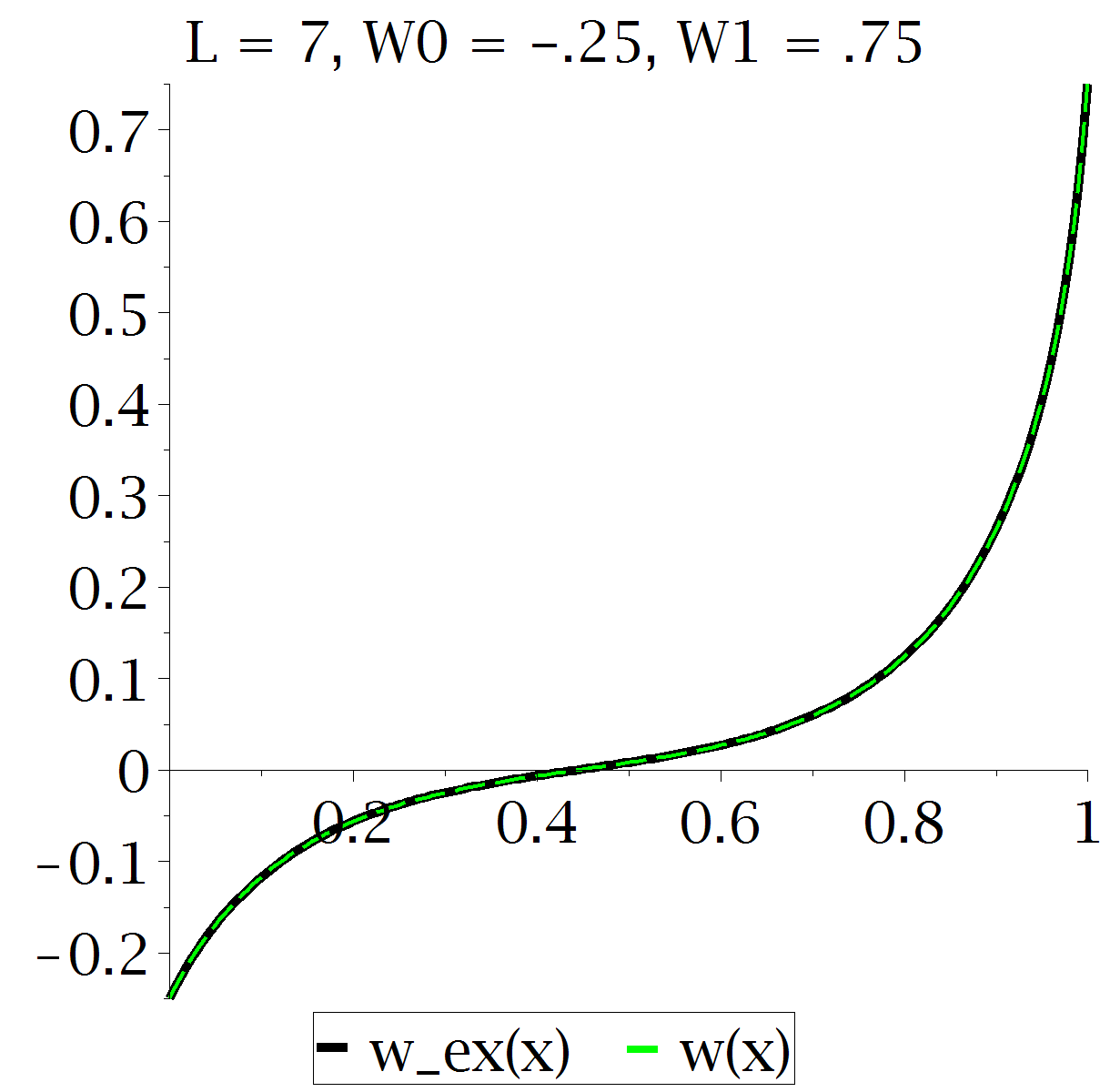}
\end{minipage}
\begin{minipage}[b]{.25\textwidth}
\includegraphics[height=.9\textwidth]{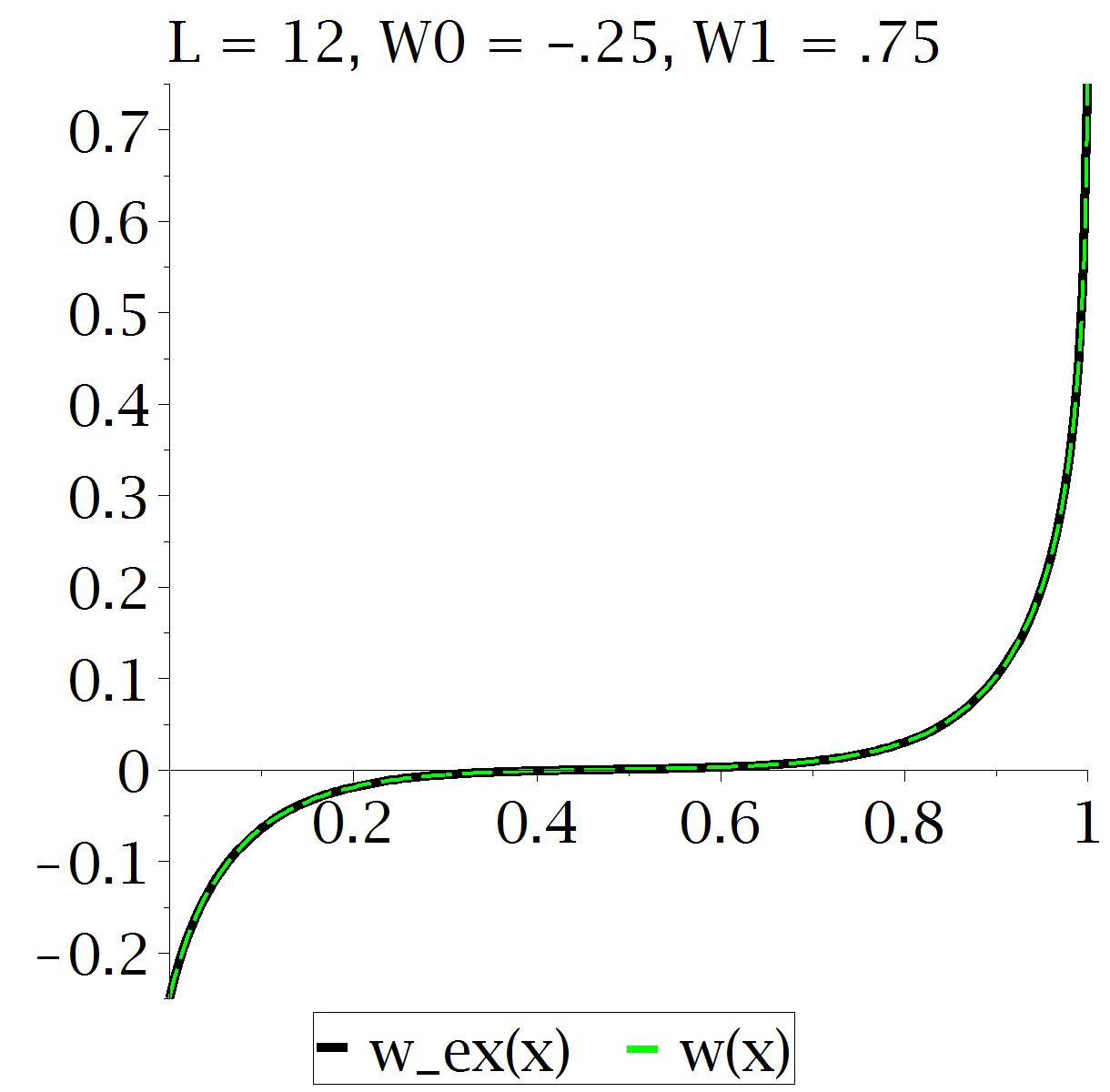}
\end{minipage}
\begin{minipage}[b]{.25\textwidth}
\includegraphics[height=.9\textwidth]{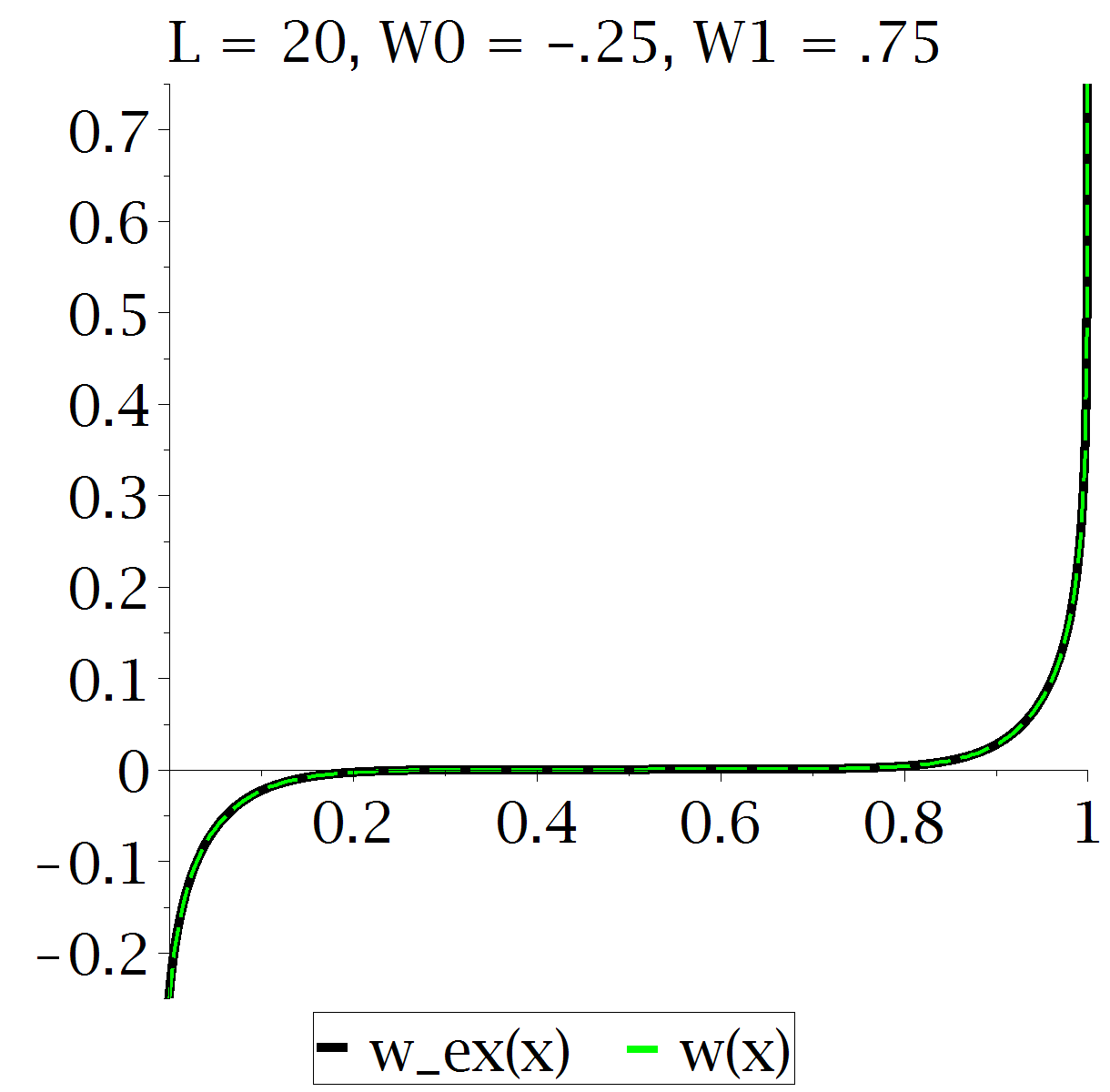}
\end{minipage}
\captionof{figure}{graph of $\ww$, for increasing values of $\LL$.}\label{fig:graph-w}
\end{minipage}

\nit We see how, as $\LL$ increases, the solution acquires boundary layers at both ends. Also, the graphs of the approximate analytical solution ---to be constructed--- overlap the exact numerical solution, their difference is very small (see Figure~\ref{fig:graph-diff}). 

Now we start analysing~\eqref{eq:gentrs-v2}. The exact solution $\ww_\exc$ has the following properties. 

\begin{enumerate}[leftmargin=5ex]
\item[(i)] 
It's strictly increasing, vanishes at exactly one point $x_0\in(0,1)$. It is negative, concave on $(0,x_0)$ and positive, convex on $(x_0,1)$.
\item[(ii)]
The restrictions of $\ww_\exc$ to $[0,x_0]$ and $[x_0,1]$ solve the classical  Troesch problem with modified boundary value. 
\item[(iii)] 
It satisfies $\ww'=\sqrt{2(\cosh(\LL\ww)-\CC)}$, with (unknown) $\CC$ less than $1$ (let $x=x_0$). For  $\xi_1:=\ww_\exc^{-1}(0.9\WW_1), \xi_0:=\ww_\exc^{-1}(0.9\WW_0)$, one has 
$
\ww_\exc'(\xi_0)\approx\re^{-0.45\LL\WW_0}\;\text{and}\;\ww_\exc'(\xi_1)\approx\re^{0.45\LL\WW_1}.
$ 
\item[] 
Hence $\ww_\exc$ has boundary layers at both $x=0$ and $x=1$, as soon as $\LL\WW_1, \LL|\WW_0|>10$. (The tangents to $\ww_\exc$ at $\xi_1, \xi_0$ make an angle less than $1^\circ$ with the vertical.)
\item[(iv)] 
The steep fall (about $x=1$) and rise (about $x=0$), makes $\ww_\exc$ almost zero on an interval containing $x_0$; the width of this interval increases with $\LL$. The graph of $\ww_\exc$ is almost flat around $x_0$, as  $\ww_\exc''(x_0)=\LL\sinh(\LL\ww_\exc(x_0))=0$. Graphically, the `almost zero in-the-middle' phenomenon can be seen in Figure~\ref{fig:graph-w}. 
\end{enumerate}

The last observation is essential because it leads to a `miracle superposition' principle for the highly non-linear equation~\eqref{eq:gentrs-v2}. The work in \S\ref{sct:ode} and~\S\ref{sct:back} allows constructing approximate solutions of Troesch-type problems (with modified right-boundary value), leaving the issue that the vanishing point $x_0$ is unknown. 

Ideally, one would extend $(\ww_\exc|_{[x_0,1]})(x)$ by $0$ on $[0,x_0]$; similarly, extend $-(\ww_\exc|_{[0,x_0]})(1-x)$ by $0$ on $[1-x_0,1]$. Instead, we approximate these by the solutions of the Troesch-type problems below (denoted $\vv_1, \vv_0$, respectively): 
\begin{m-eqn}{
\vv''=\LL\cdot\sinh(\LL\vv),\;\,\vv(0)=0,\;\vv(1)=\WW_1\;\; \text{(resp. $-\WW_0$)}.
}\label{eq:dec}
\end{m-eqn} 
They are both defined on $[0,1]$ and possess boundary layer for sufficiently large $\LL$; as explained, numerically, this means that $\LL\WW_1, \LL\WW_0>10$. 


\subsubsection{`miracle superposition'} 

The boundary layer causes that $\vv_1, \vv_0$ almost vanish away from $x=1$, allowing to treat them as approximations of the extensions-by-zero above.
\\[1ex] \begin{minipage}[c]{.97\textwidth}\centering 
\begin{minipage}[c]{.60\textwidth} 
We clarify the discussion with a picture: we plotted the graphs of $\vv_1, \vv_0$ (actually of their analytical approximations). One can see that the addition yields the graph in the middle of Figure~\ref{fig:graph-w}. 
\end{minipage}\hskip3ex
\begin{minipage}[c]{.30\textwidth}
\includegraphics[height=.80\textwidth]{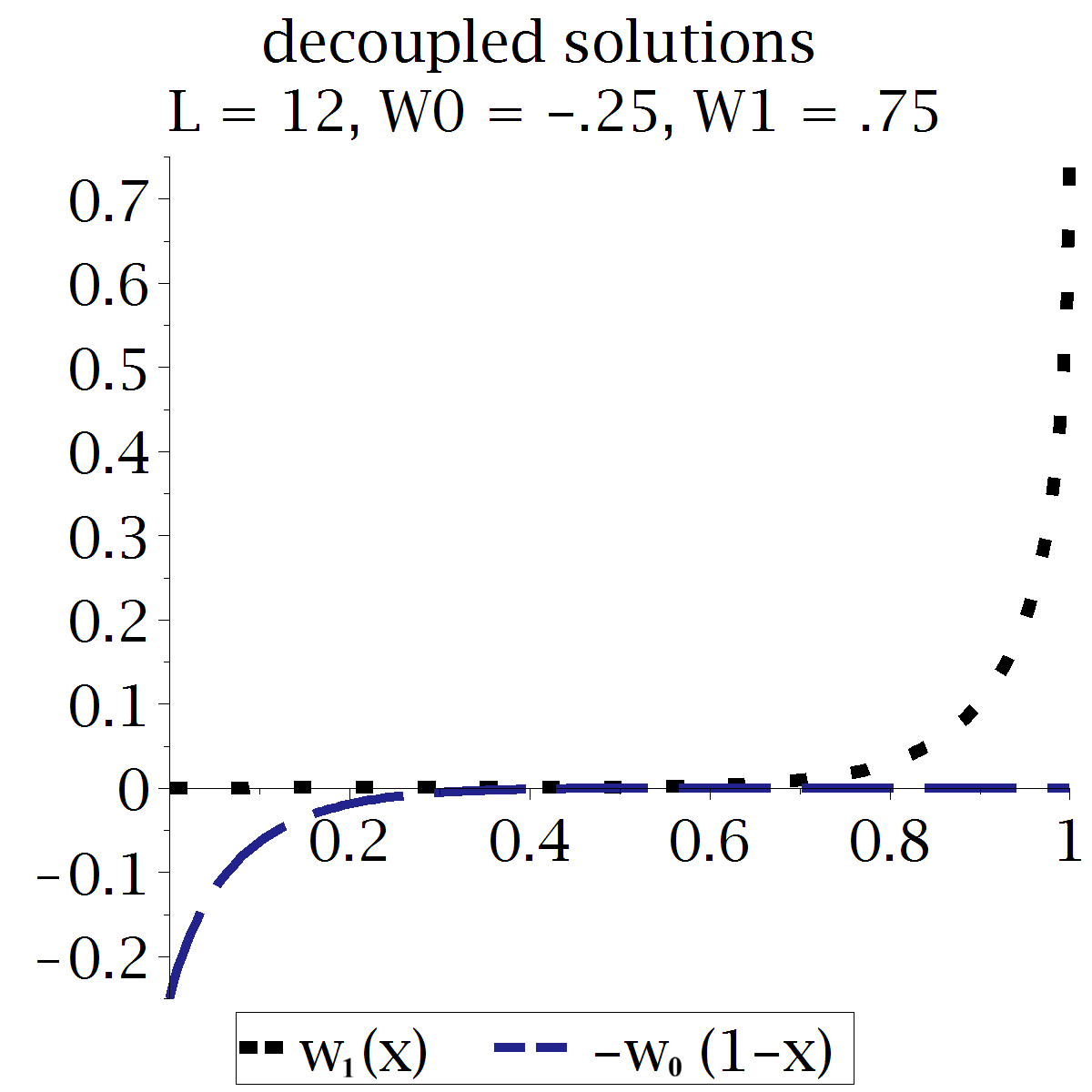}
\end{minipage}
\end{minipage}

\nit Explicit analytic estimates of $\{x\mid |\vv_1|, |\vv_0|<\veps \}$ are tedious; see~\cite[\S2]{csfx}) how to squeeze $\vv_1, \vv_0$ between explicit upper/lower bounds. Nevertheless, the guiding idea of the superposition is straightforward: for large $\LL$, $\vv_1, \vv_0$ have boundary layer near $x=1$ and almost vanish for $x\in[0,0.5]$. 
Thus, for $x>0.5$, we have $\vv_0(1-x)\approx0$ and $\vv_0''(1-x)=\sinh(\LL\vv_0(1-x))\approx0$: 
$$\LL\sinh(\LL\vv(x))\approx \LL\sinh(\LL\vv_1(x))=\vv_1''(x)\approx\vv_1''(x)-\vv_0''(1-x)=\vv''(x).$$
Similar reasoning shows that this approximation is true for $x<0.5$, too.

\begin{m-proposition}\label{prop:super}
The function $\vv(x):=\vv_1(x)-\vv_0(1-x)$ is an approximate solution of the {\bvp}~\eqref{eq:gentrs-v2} and the accuracy of the estimate improves for increasing $\LL$.

Consequently, let $\ww_1, \ww_0$ be the explicit approximate solutions to~\eqref{eq:dec}, constructed in \S\ref{ssct:dec}. Then $\ww(x):=\ww_1(x)-\ww_0(1-x)$ is an approximate solution to~\eqref{eq:gentrs-v2}. 
\end{m-proposition}

\begin{m-proof}
As soon as $\vv''-\LL\sinh(\LL\vv)$ is small (in $\max$-norm), the continuos dependence of solutions on parameters implies that $\vv$ is close to the exact solution $\ww_\exc$.
\end{m-proof}


\subsection{Solving the decoupled {\bvp}s}\label{ssct:dec} 

We reduced the problem of approximating the solution of~\eqref{eq:gentrs-v2} to solving the `decoupled' pair of {\bvp}s in Proposition~\ref{prop:super}. This is where our `analytical toolbox' is essential. 

\begin{m-proposition}\label{prop:super-w}
The functions 
$$\begin{array}{l}
\ww_1(x)=\frac{2}{\LL}\ln\bigg[
\CC_1\cdot \frac{\re^{\frac{\LL\WW_1}{2}}+\CC_1\cdot\tanh\big[ \frac{\LL\CC_1(1-x)}{2} \big]}
{\CC_1+\re^{\frac{\LL\WW_1}{2}}\cdot\tanh\big[ \frac{\LL\CC_1(1-x)}{2} \big]}
\bigg],\,
\CC_1=\tanh
\frac{\tanh^{-1}\big(\re^{\frac{-\LL\WW_1}{2}}\big)+\frac{\LL(\sinh\LL-\LL)}{4\cosh^2(\LL/2)}}{1-\frac{\LL}{2\cosh^2(\LL/2)}},
\\[2ex] 
\ww_0(x)=\frac{2}{\LL}\ln\bigg[
\CC_0\cdot \frac{\re^{\frac{-\LL\WW_0}{2}}+\CC_0\cdot\tanh\big[ \frac{\LL\CC_0(1-x)}{2} \big]}
{\CC_0+\re^{\frac{-\LL\WW_0}{2}}\cdot\tanh\Big[ \frac{\LL\CC_0(1-x)}{2} \big]}
\bigg],\,
\CC_0=\tanh
\frac{\tanh^{-1}\big(\re^{\frac{\LL\WW_0}{2}}\big)+\frac{\LL(\sinh\LL-\LL)}{4\cosh^2(\LL/2)}}{1-\frac{\LL}{2\cosh^2(\LL/2)}},
\end{array}
$$ 
are both upper envelopes of the exact solutions of~\eqref{eq:dec}. In fact, we have 
$$\ww_\exc(x)\les\ww_1(x)\;\;\text{and}\;\;-\ww_\exc(1-x)\les-\ww_0(1-x),\;\;\text{for}\;x\in[0,1].$$
We declare that the function $\;\ww(x)=\ww_1(x)-\ww_0(1-x)\;$ 
is an approximate solution of~\eqref{eq:gentrs-v2}. It (always) satisfies $\ww(0)>\WW_0,\; \ww(1)<\WW_1$.  
\end{m-proposition}

\begin{m-proof}
We analyse the case of $\WW_1$, the other is analogous. By integrating the defining {\ode} we obtain 
$\ww'=\sqrt{\ee^{\LL\ww}+\ee^{-\LL\ww}-2\CC^2}.$ Now we apply the algorithm described in \S\ref{sct:ode}.  The term $\re^{\LL\ww}$ is responsible for the large derivative ---equivalently, the large curvature $\yy''(1)$---, for $\LL\gg0$. This imposes the change of coordinates $z=\re^{\LL\yy/2}$, which yields:
$$
\zz'=(\LL/2)\cdot\sqrt{\zz^4+1-2\CC^2\zz^2},\;\zz(0)=1,\,\zz(1)=\re^{\frac{\LL\WW_1}{2}}.
$$ 
The steep variation is caused by $\zz(1)\gg0$, so we need a quadratic large-$\zz$ expansion:  
 $\sqrt{\zz^4+1-2\CC^2\zz^2}=\zz^2\sqrt{1-(2\CC^2\zz^2-\zz^{-4})}\approx\zz^2(1-\CC^2\zz^{-2})$. 

The estimate is a `$\ges$'-inequality, so 
$\;\zz'=(\LL/2)\cdot(\zz-\CC)(\zz+\CC),\;\zz(0)=1, \zz(1)=\re^{\frac{\LL\WW_1}{2}}\;$ 
yields the upper envelope $\ww_1\ges\ww_\exc$: both have the same initial value, but the slope of the exact solution is greater. The explicit form of $\ww_1$ is given in the Proposition, and~\eqref{eq:cc5} yields the approximate value $\CC=\CC_1$. 
Finally, we compute $\ww(1)=\ww_1(1)-\ww_0(0)=\WW_1-\ww_0(0)<\WW_1$; for the last step, we used $\ww_0(0)>0$, as $\ww_0$ is an upper bound for $\ww_\exc$.
\end{m-proof}


\subsection{Numerical tests}

We verify that the approximate solution $\ww$ is (highly) precise and it improves with $\LL$ (the difference $\ww-\ww_\exc$ decreases). Recall that the error $O(\exp(-\LL))$ (see \cite[Prop. 2.6]{csfx}) of the approximate solutions $\ww_1, \ww_0$ of the decoupled {\bvp}s was confronted against numerical data. Figure~\ref{fig:graph-w} shows that, for our generalized equation, the graph of $\ww$ overlaps the (numerical) solution $\ww_\exc$. We distinguish three methods to estimate their difference. 

\subsubsection{solving the \bvp} 
Typically, this works for low values of $\LL\cdot\max(|\WW_0|,\WW_1)$. 
\\[1ex]\begin{minipage}{.99\textwidth}\centering 
\begin{minipage}[c]{.455\textwidth} 
The double boundary layer makes~\eqref{eq:gentrs-v2} difficult. The classical Troesch equation ---it possesses only one boundary layer--- is notoriously known for being extremely sensitive. Thus our generalized problem is (about) twice so.
\end{minipage}\qquad
\begin{minipage}[c]{.24\textwidth}
\includegraphics[height=.9\textwidth]{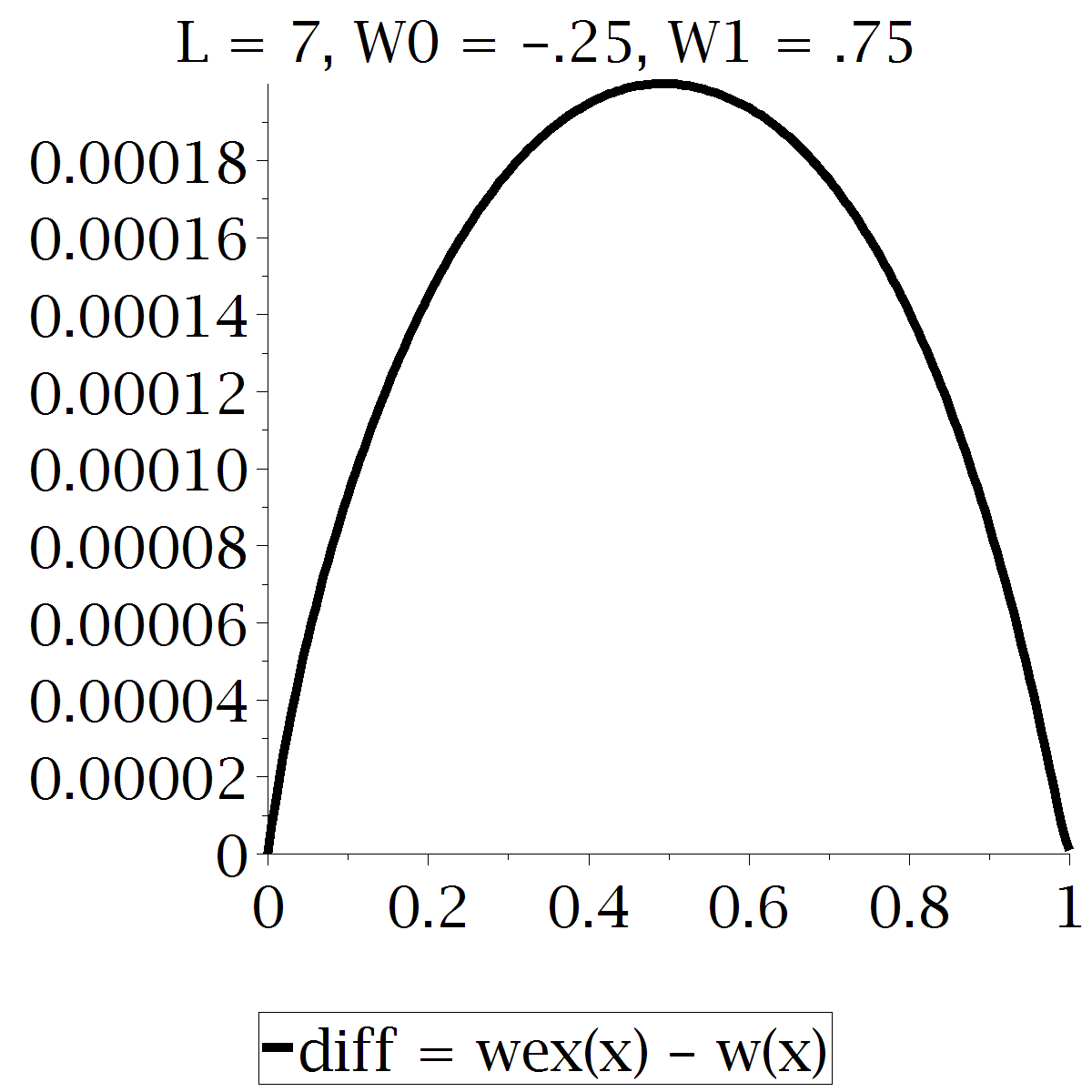}
\end{minipage}
\begin{minipage}[c]{.24\textwidth}
\includegraphics[height=.9\textwidth]{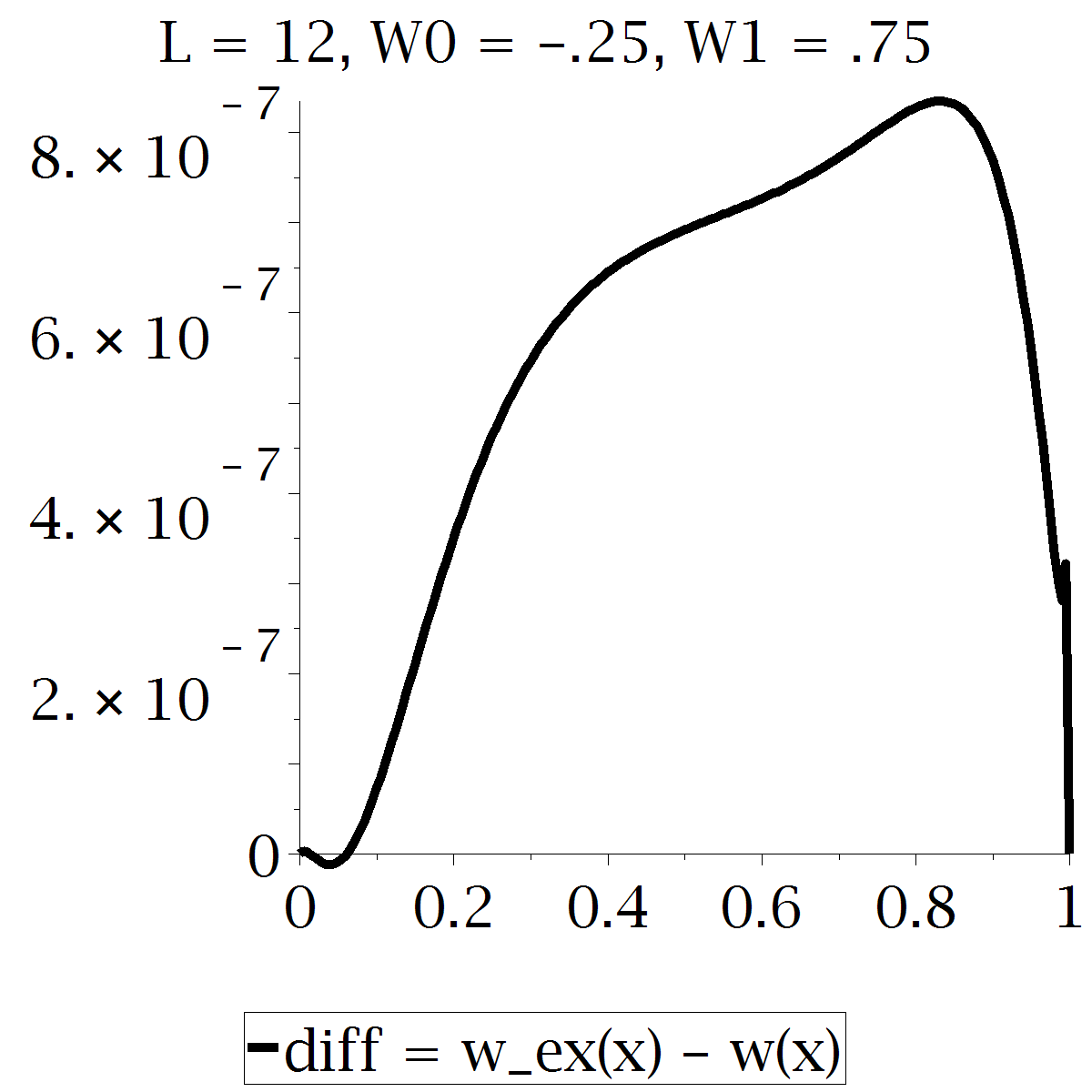}
\end{minipage}
\captionof{figure}{estimating the error $\ww-\ww_\exc$, for increasing values of $\LL$.}\label{fig:graph-diff}
\end{minipage}


\subsubsection{reducing the order} 

The values of $\CC_0, \CC_1$ in Proposition~\ref{prop:super-w} are very accurate, especially in the boundary layer case, allowing to transform the {\bvp}~\eqref{eq:gentrs-v2} into the {\ivp}s: 
\begin{m-eqn}{
\begin{array}{l}
\zz'=(\LL/2)\cdot\sqrt{\zz^4+1-2\CC^2\zz^2},
\\ 
\zz=\exp(\LL\ww_\exc/2),
\end{array}\quad 
\bigg\{\begin{array}{rll}
\CC=\CC^{(x=1)}_{\exc},\;\zz(1)=\re^{\frac{\LL\WW_1}{2}},&\;\text{for}\;x=1;
&\leadsto\text{soln.}\;Z^{(1)}
\\ 
\CC=\CC^{(x=0)}_{\exc},\;\zz(0)=\re^{\frac{\LL\WW_0}{2}},&\;\text{for}\;x=0;
&\leadsto\text{soln.}\;Z^{(0)}.
\end{array}
}\label{eq:gentrs-zz}
\end{m-eqn}
The values $\CC^{(x=1)}_{\exc}, \CC^{(x=0)}_{\exc}$ are obtained by `shooting', starting with $\CC_1, \CC_0$, respectively. The {\ivp}s are much easier to solve numerically, for large parameters, too. Thus can be compared with $\ww$ in Proposition~\ref{prop:super-w}. The bottom right cell of Table~\ref{tab:c1c0} illustrates the sensitivity of the {\ivp} on the data. To improve the accuracy, one should continue the shooting. 

\begin{center}\renewcommand{\arraystretch}{1.1}
\textscale{.85}{\begin{tabular}[H]{|c|c|c|c|c|}
\hline
$(\LL,\WW_0,\WW_1)$&(20,-0.25,0.75)&(40,-0.25,0.75)&(15, -0.5, 2)&(30, -0.5, 2)
\\\hline  
$\CC^{(x=1)}_\exc$&$(1-2.3826.10^{-8})\CC_1$&$(1-5.8838.10^{-17})\CC_1$&$(1-4.0577.10^{-6})\CC_1$&$(1-1.3084.10^{-12})\CC_1$
\\ 
$\frac{2}{\LL}\ln\big(Z^{(1)}(0)\big)-\WW_0$&$8.8\cdot10^{-6}$&$3.9\cdot10^{-5}$&$5.2\cdot10^{-5}$&$5.5\cdot10^{-4}$
\\\hline  
$\CC^{(x=0)}_\exc$&$(1-2.4447.10^{-8})\CC_0$ &$(1-5.8804.10^{-17})\CC_0$&$(1-4.0859.10^{-6})\CC_0$&$(1-1.308614.10^{-12})\CC_0$
\\ 
$\WW_1-\frac{2}{\LL}\ln\big(Z^{(0)}(1)\big)$&$8.8\cdot10^{-4}$&$1.1\cdot10^{-1}$&$6.6\cdot10^{-2}$&$6.3\cdot10^{-1}$
\\\hline
\end{tabular}}
\captionof{table}{approximate values of $\CC_\exc$ and error at the (opposite) end}\label{tab:c1c0}
\renewcommand{\arraystretch}{1}\end{center}

\nit We approximate $\ww_\exc$ by  
$\,Z(x)=\bigg\{\begin{array}{ll}Z^{(0)}(x),&\,x\in[0,0.5];\\ Z^{(1)}(x),&\,x\in[0.5,1],\end{array}\,$ and compute $\ww(x)-Z(x)$. 
\\[1ex]\begin{minipage}{.975\textwidth}\centering
\begin{minipage}[b]{.24\textwidth} 
\includegraphics[height=.9\textwidth]{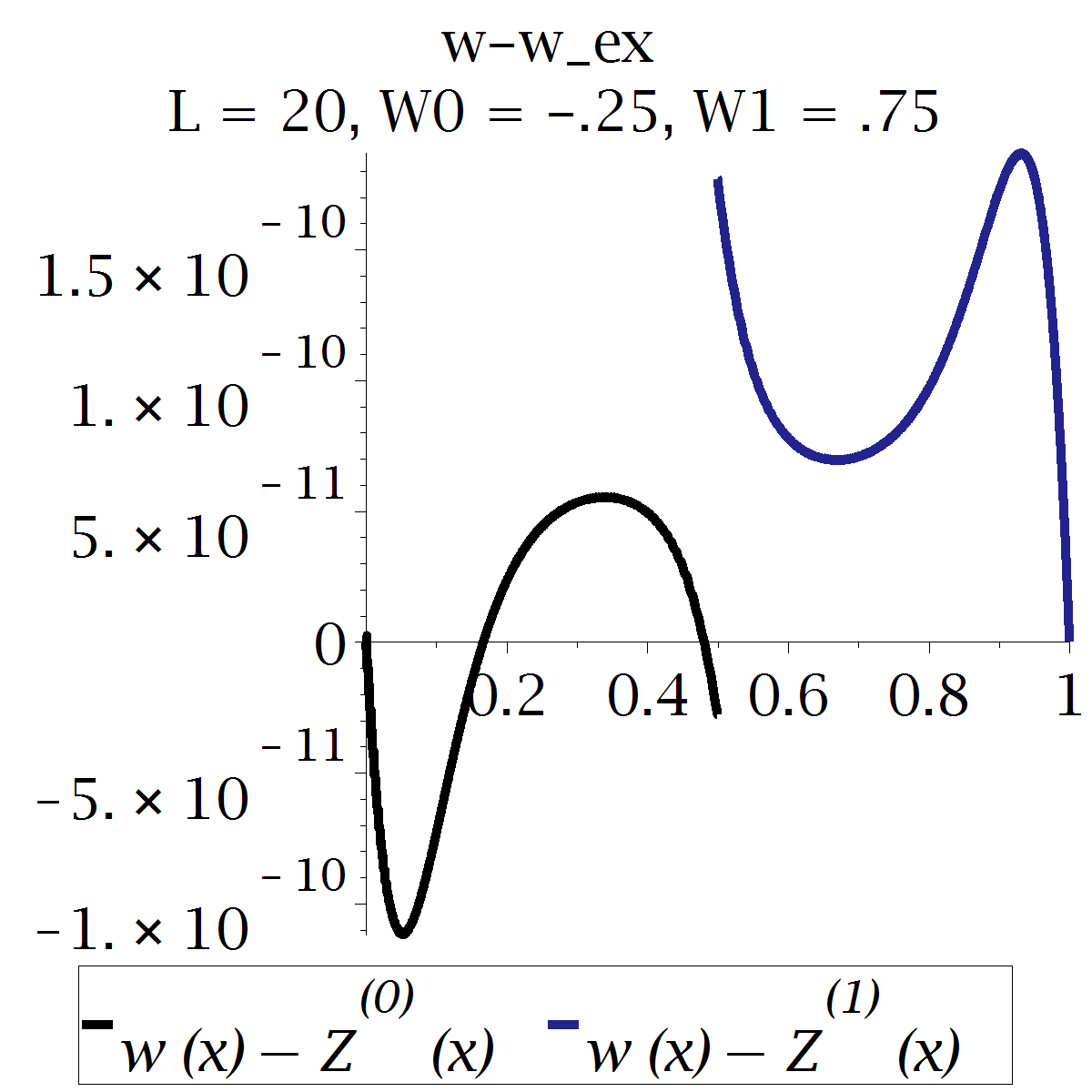}
\end{minipage}\hskip-.5ex
\begin{minipage}[b]{.24\textwidth} 
\includegraphics[height=.9\textwidth]{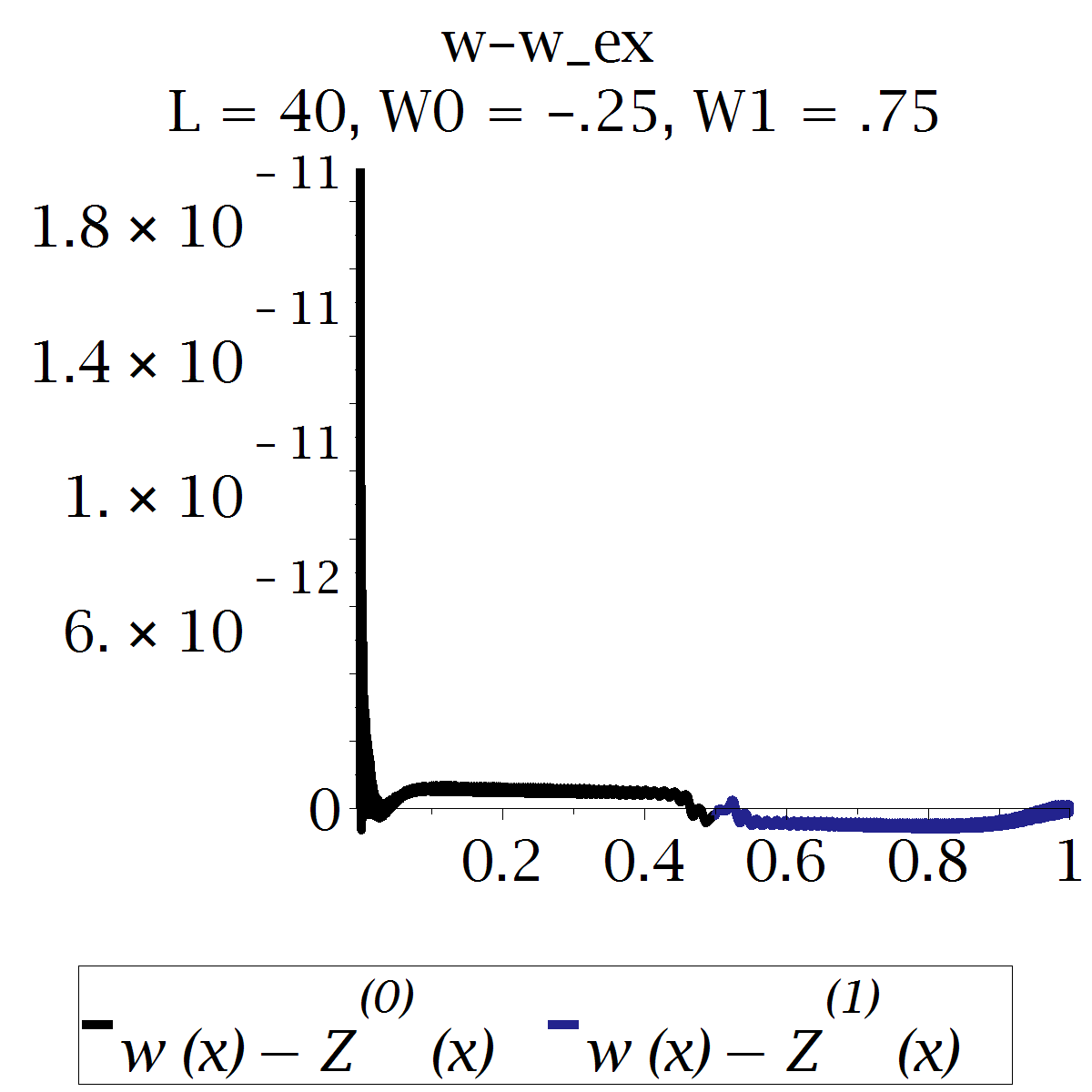}
\end{minipage}
\begin{minipage}[b]{.24\textwidth} 
\includegraphics[height=.9\textwidth]{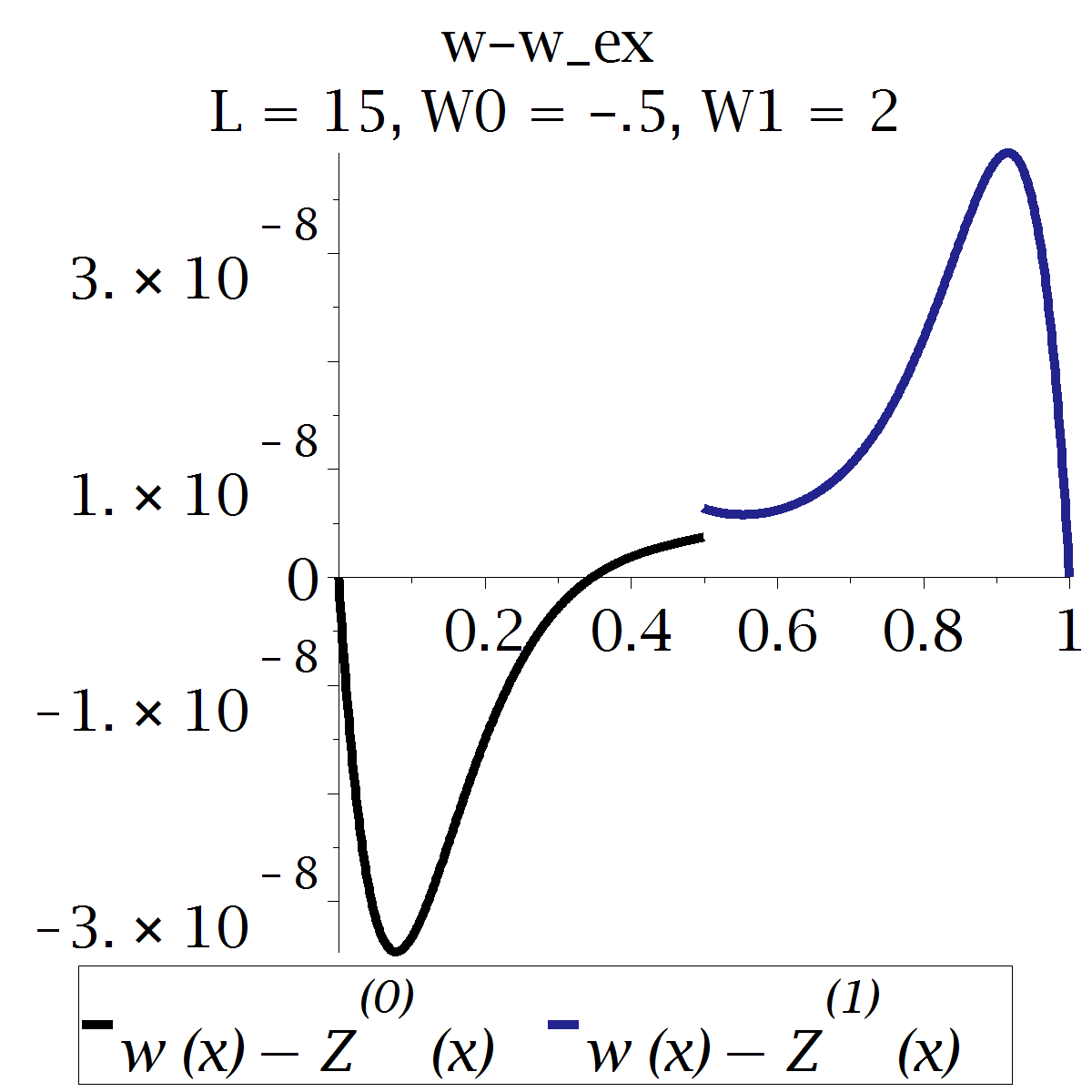}
\end{minipage}\hskip-.5ex
\begin{minipage}[b]{.24\textwidth} 
\includegraphics[height=.9\textwidth]{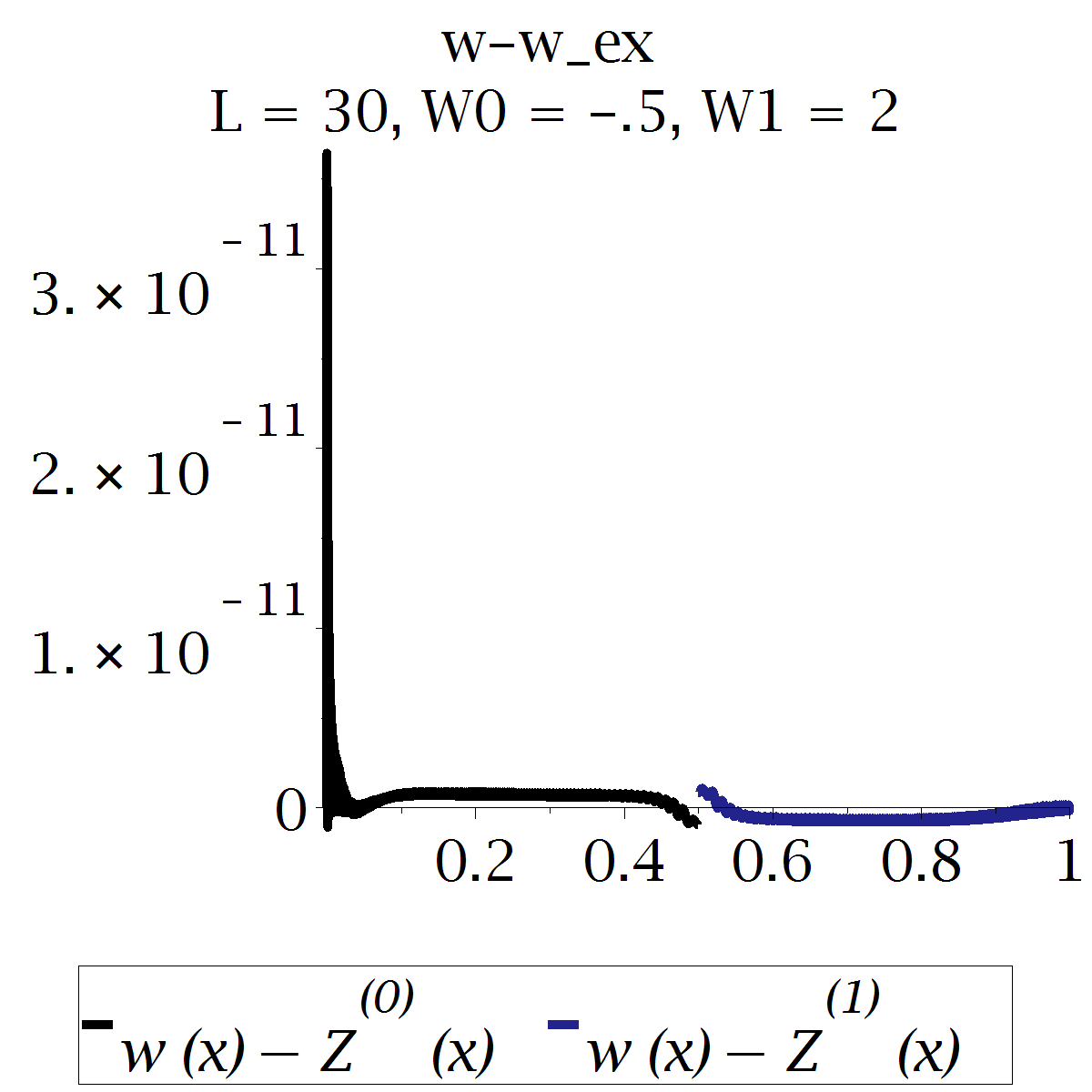}
\end{minipage}
\captionof{figure}{estimating the error $\ww(x)-\ww_\exc(x)$. }
\end{minipage}

\subsubsection{computing the residue}

No hard-to-compute numerical solutions are needed. Rather, we compute the residue ---in our situation $\resd_\ww(x):=\ww''(x)-\LL\cdot\sinh(\LL\ww(x))$--- and use the last remark in Proposition~\ref{prop:super-w}. 

When both $\ww_1, \ww_0$ have boundary layer, the computationally difficult case, the superposition principle implies  
$\resd_\ww(x)\approx\resd_{\ww_1}(x)-\resd_{\ww_0}(1-x).$ 
So the residue is (very) small away from $x=0, 1$, so $\ww$ almost solves the {\ode} $\ww''=\LL\sinh(\LL\ww)$.  (See Figure~\ref{fig:graph-res}. We included examples where the boundary values are not of the form~\eqref{eq:Wp}). 
\\[1ex]\begin{minipage}{.99\textwidth}\centering 
\begin{minipage}[b]{.19\textwidth}
\includegraphics[height=.95\textwidth]{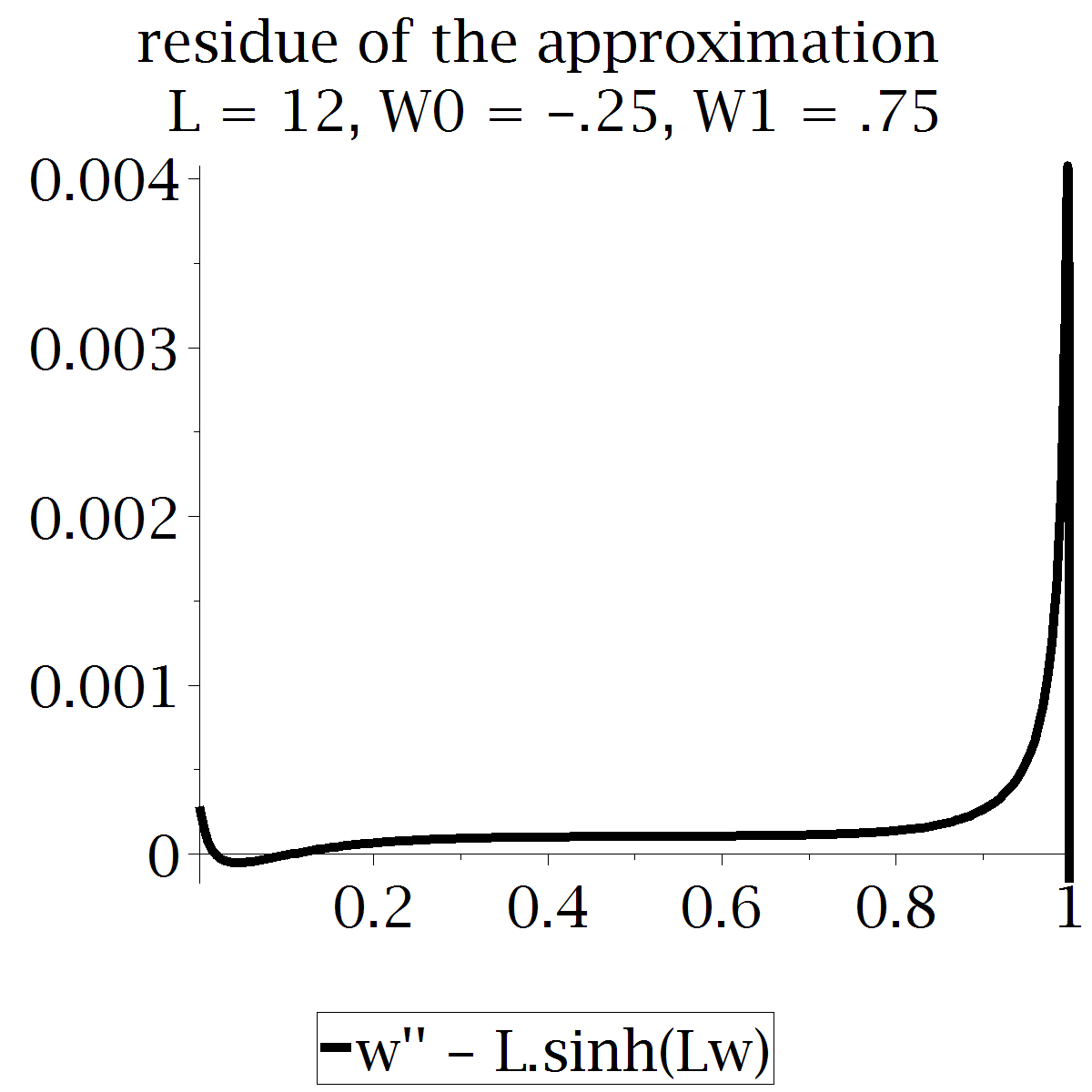}
\end{minipage}
\begin{minipage}[b]{.19\textwidth}
\includegraphics[height=.95\textwidth]{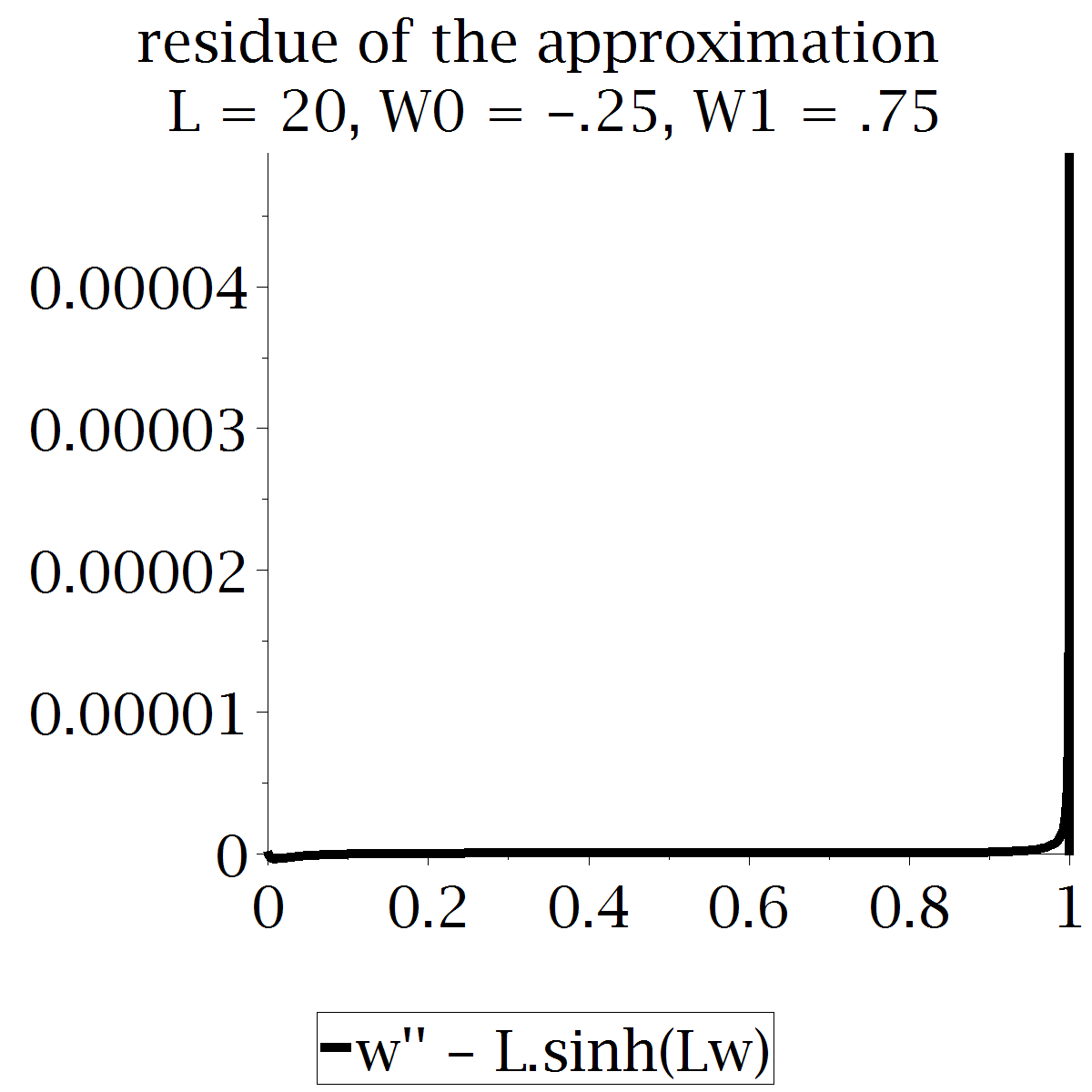}
\end{minipage}
\begin{minipage}[b]{.19\textwidth}
\includegraphics[height=.95\textwidth]{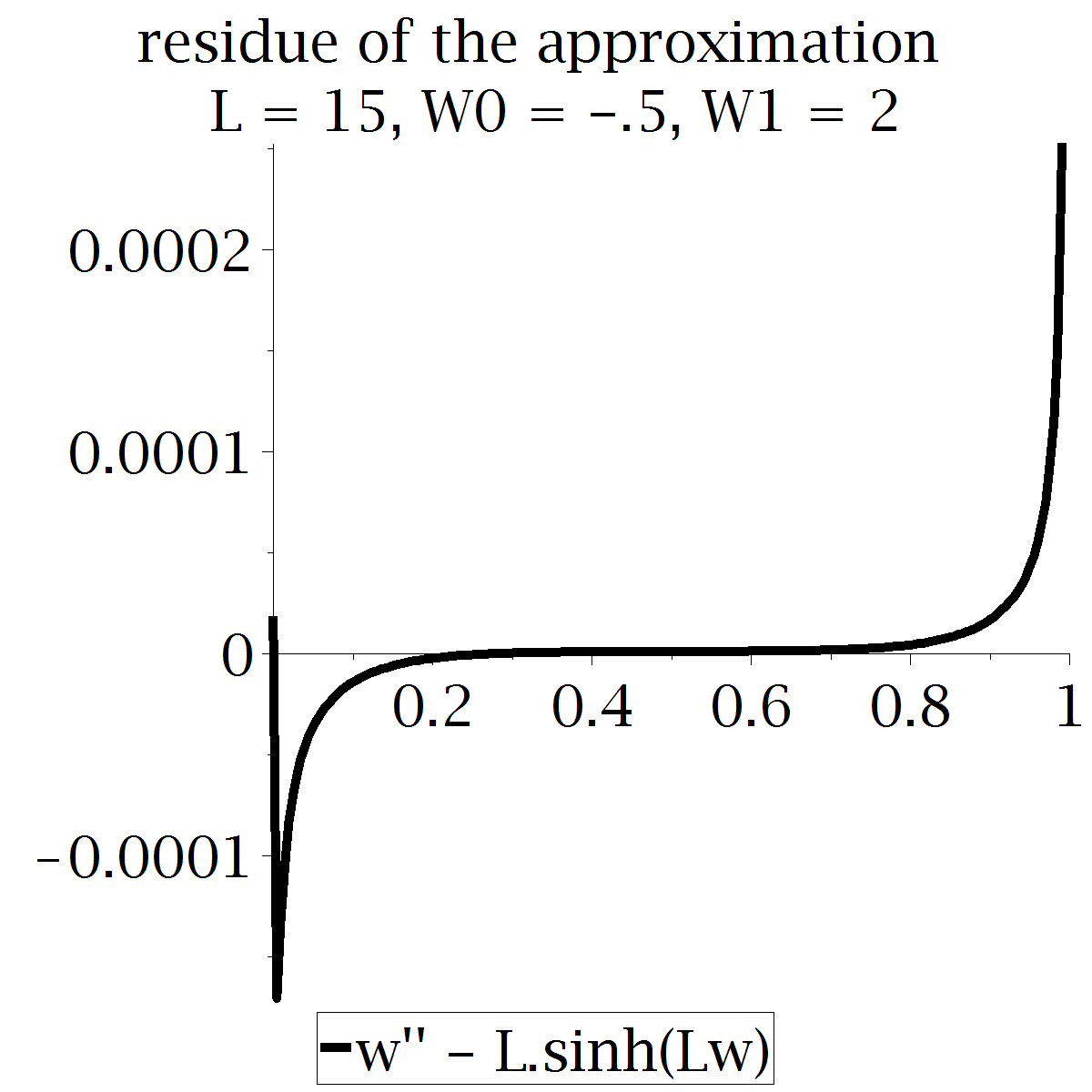}
\end{minipage}
\begin{minipage}[b]{.19\textwidth}
\includegraphics[width=.95\textwidth]{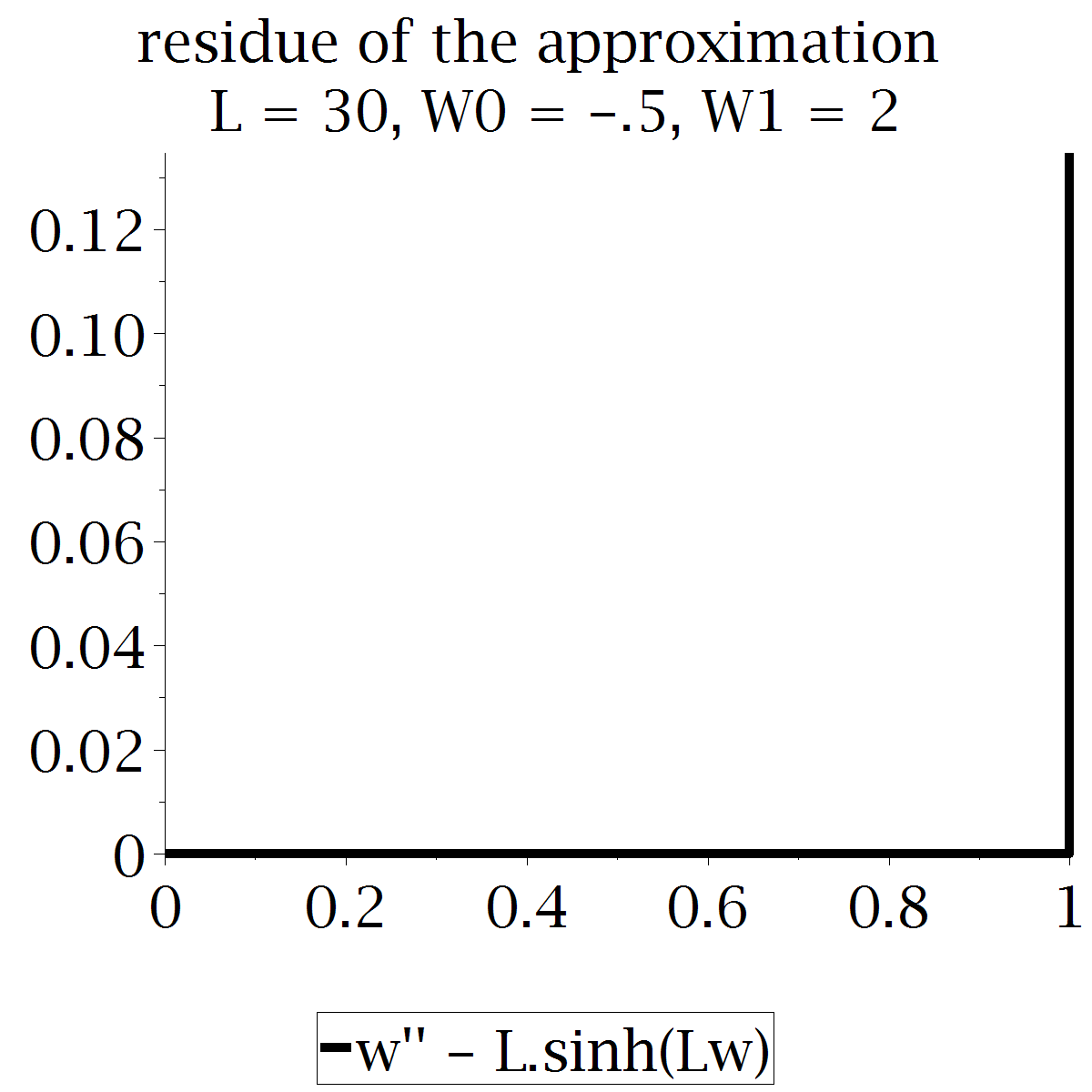}
\end{minipage}
\begin{minipage}[b]{.19\textwidth}
\includegraphics[width=.95\textwidth]{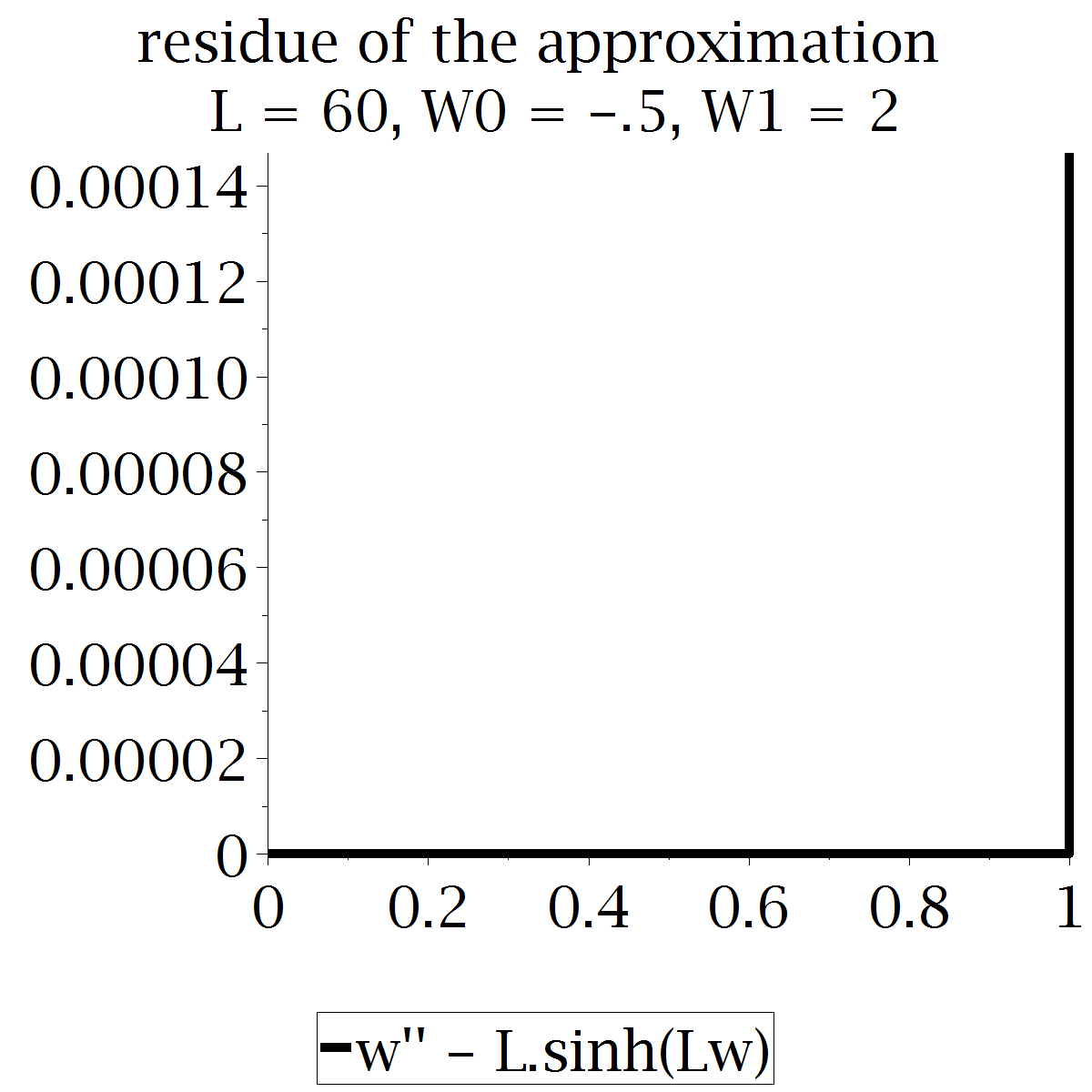}
\end{minipage}
\captionof{figure}{graphs of the residue function.}\label{fig:graph-res}
\end{minipage}

\nit To conclude that $\ww$ approximates well $\ww_\exc$, we need to verify that $\ww-\ww_\exc$ is small near $x=1$ and $x=0$ (where, unfortunately, the residue increases). This is where we need the inequalities  and $\ww(0)>\WW_0=\ww_\exc(0)$. Since $\ww_1$ is a certain (global) upper envelope of $\ww_\exc$, the inequality $\ww(1)<\WW_1=\ww_\exc(1)$ implies that, in a (very) small neighbourhood of $x=1$, we have 
$\;0<\ww_\exc(x)-\ww(x)<\ww_1(x)-\ww(x)=\ww_0(1-x);\;$  
we know that the latter is very small. Similar argument works near the $x=0$ end.
\\[1ex]\begin{minipage}{.99\textwidth}\centering 
\begin{minipage}[b]{.19\textwidth}
\includegraphics[height=.95\textwidth]{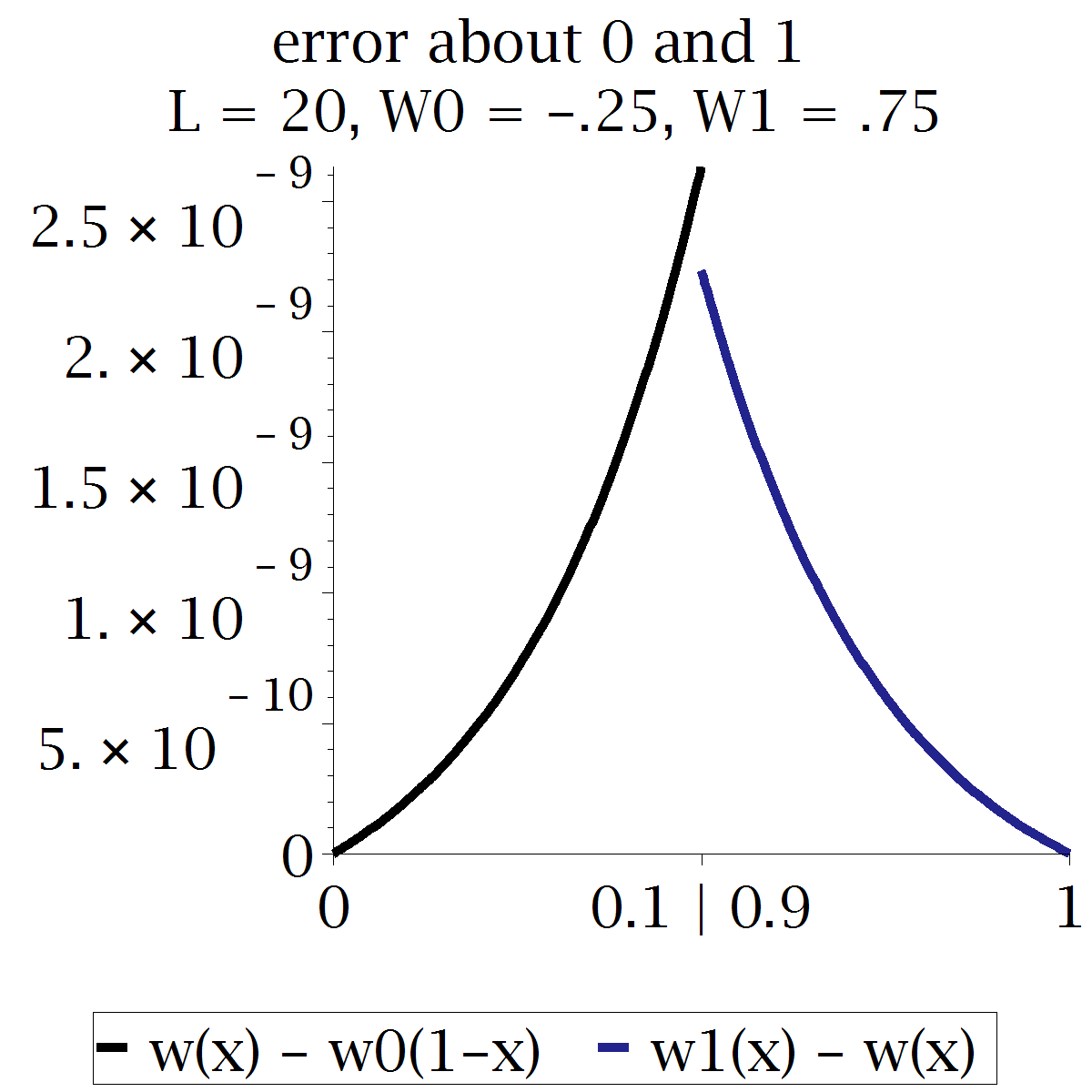}
\end{minipage}
\begin{minipage}[b]{.19\textwidth}
\includegraphics[height=.95\textwidth]{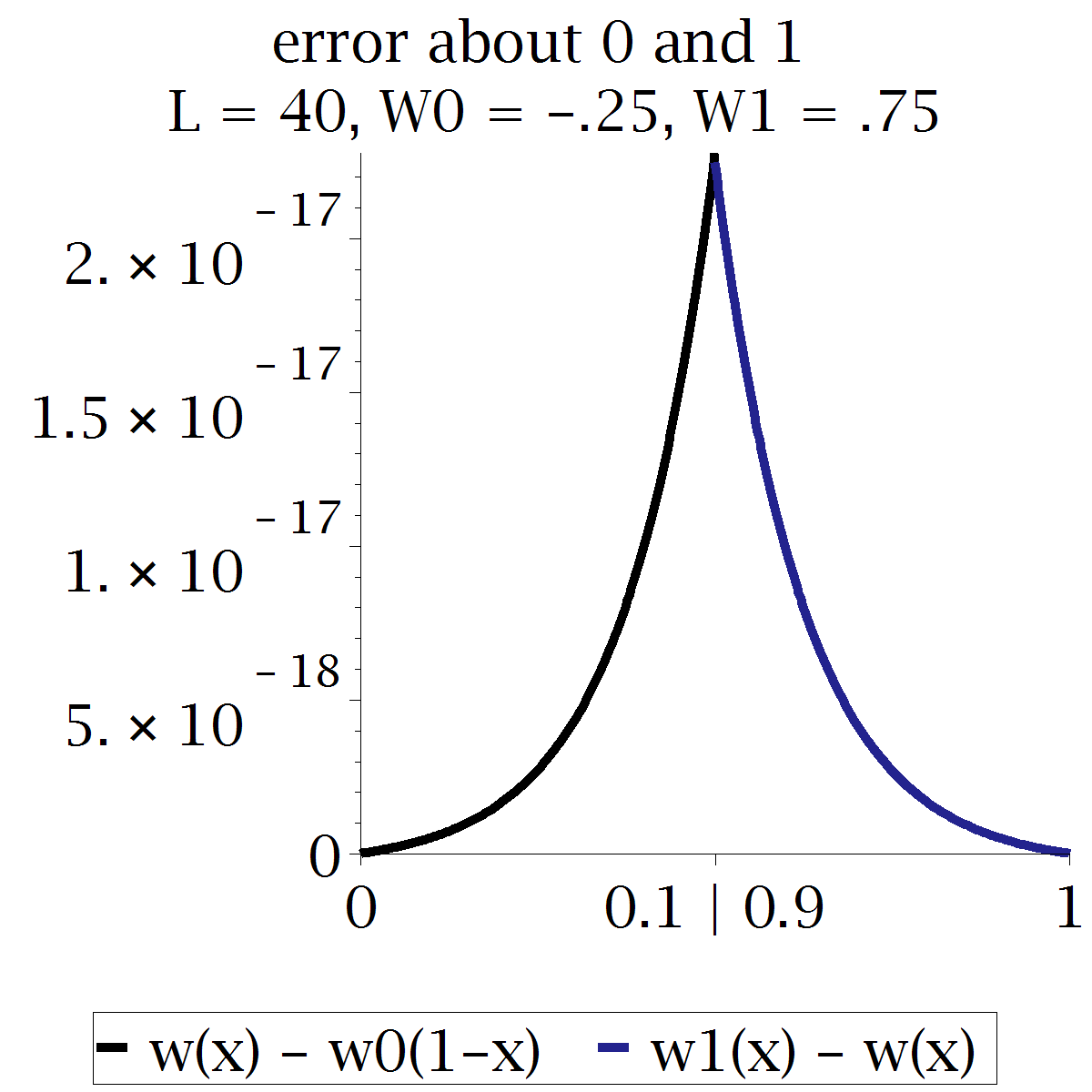}
\end{minipage}
\begin{minipage}[b]{.19\textwidth}
\includegraphics[height=.95\textwidth]{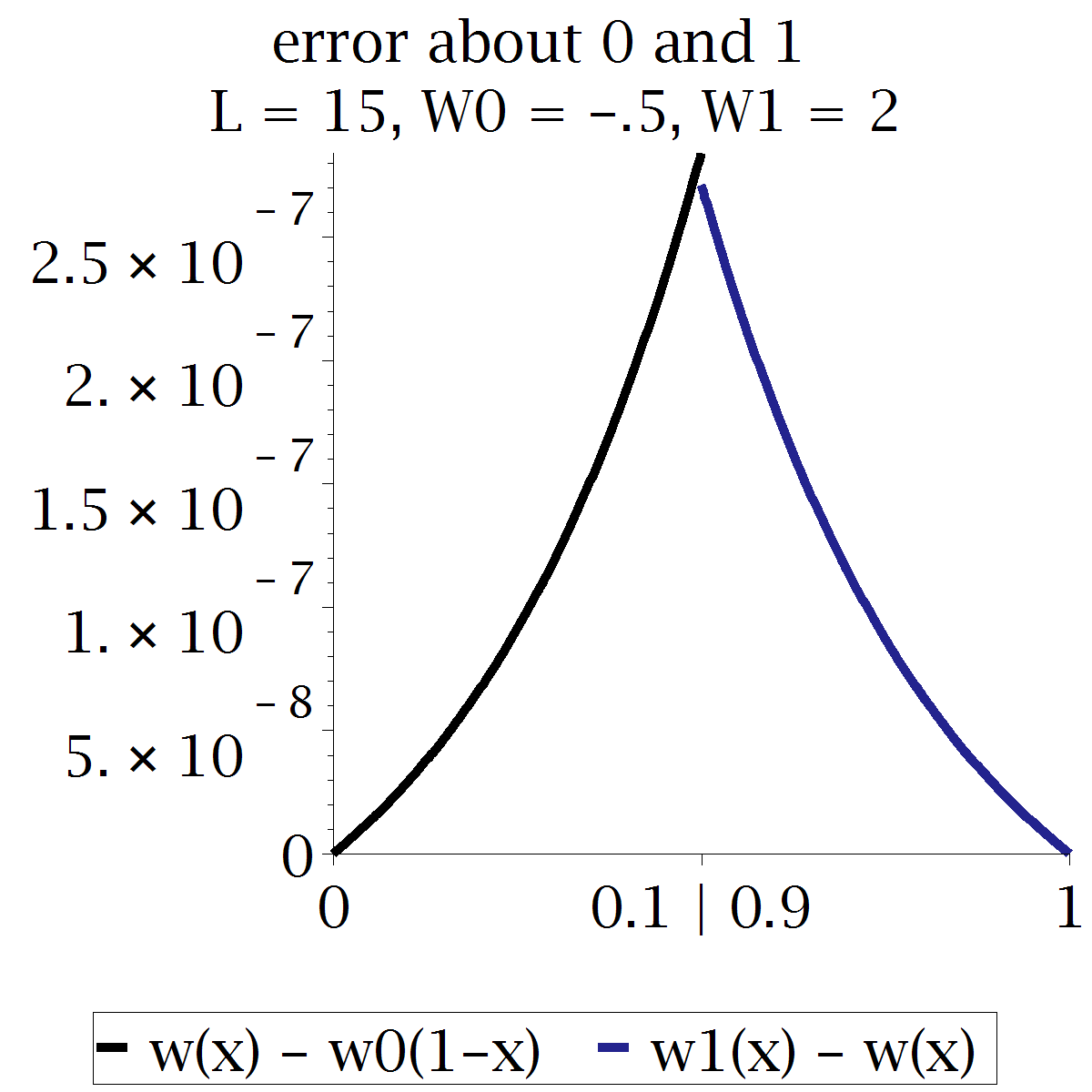}
\end{minipage}
\begin{minipage}[b]{.19\textwidth}
\includegraphics[width=.95\textwidth]{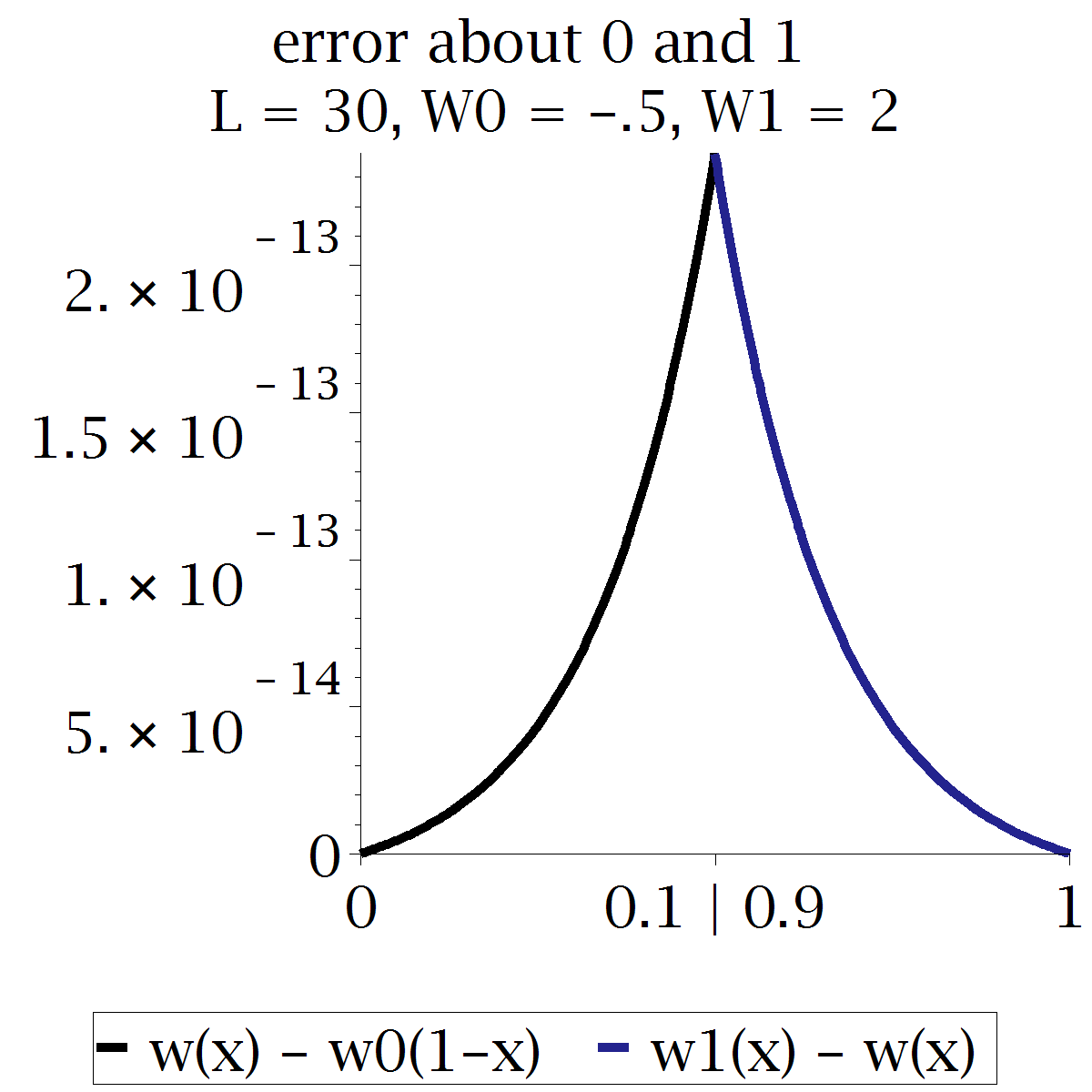}
\end{minipage}
\begin{minipage}[b]{.19\textwidth}
\includegraphics[width=.95\textwidth]{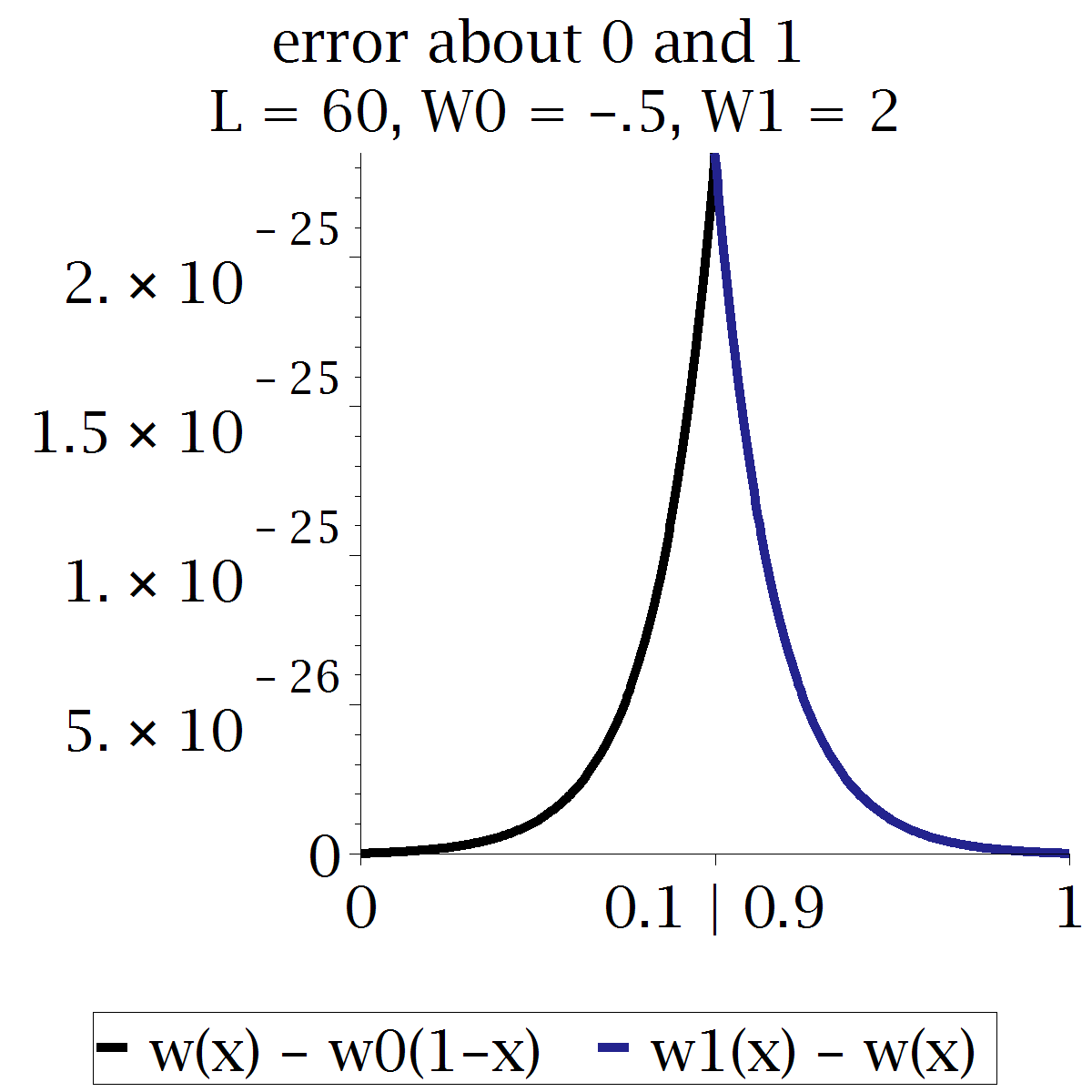}
\end{minipage}
\captionof{figure}{estimating the error about $x=0$ and $x=1$.}\label{fig:graph-diffend}
\end{minipage}

\nit We observe that the order of precision (exponents of $10$) are consistent with those in Figure~\ref{fig:graph-diff}. This approach is probably the most effective to verify the accuracy of the estimate $\ww\approx\ww_\exc$.


\subsection{Back to the roots}\label{ssct:phys}

The motivation for this work was to include charge separation into Troesch's problem. The ion/electron density functions are $\;\NN_i(x)=\re^{\LL\ww(x)},\;\,\nn_\ee=\re^{-\LL\ww(x)}.$

\begin{minipage}{.99\textwidth}\centering 
\begin{minipage}[b]{.32\textwidth}
\includegraphics[height=.7\textwidth]{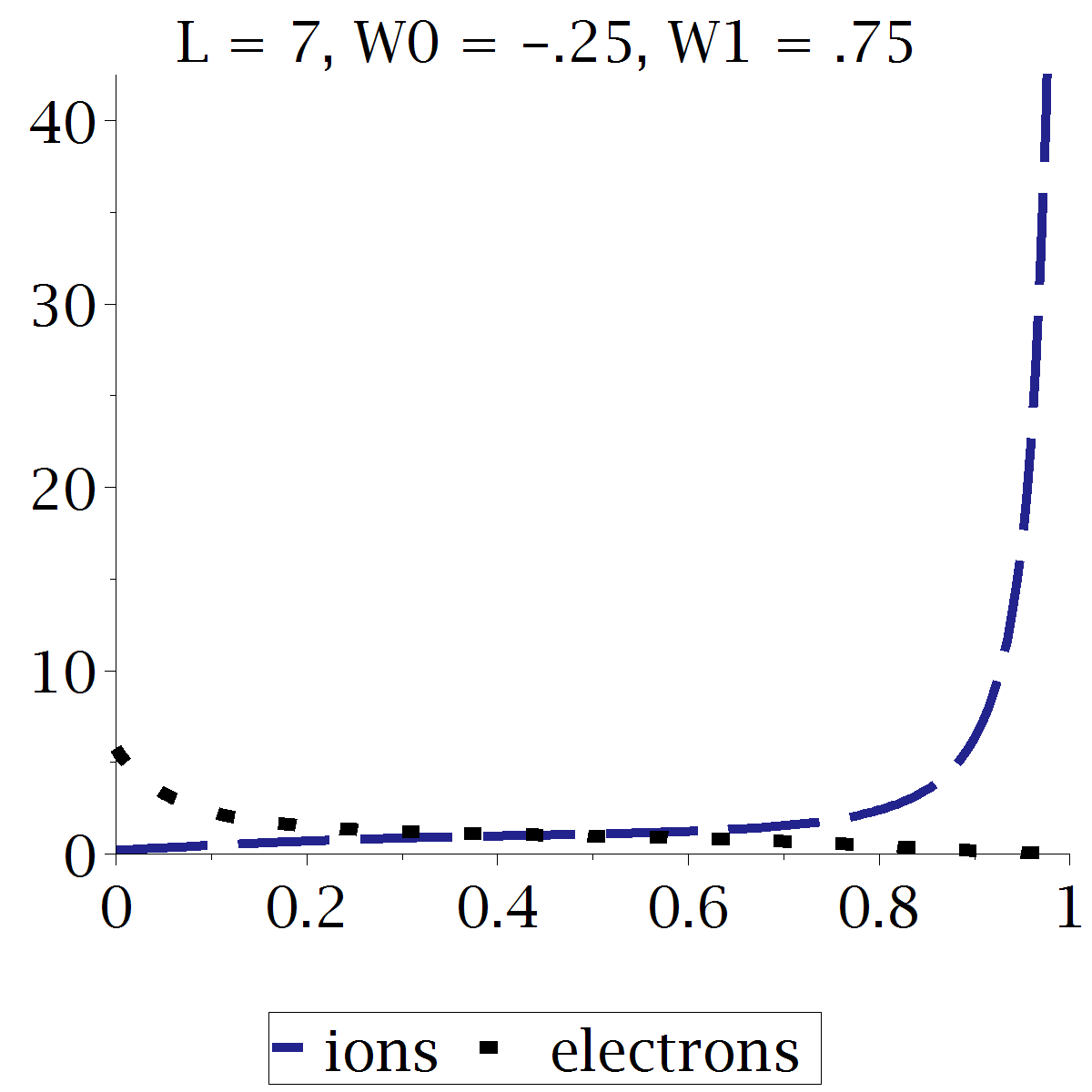}
\end{minipage}\hskip0ex 
\begin{minipage}[b]{.32\textwidth}
\includegraphics[height=.7\textwidth]{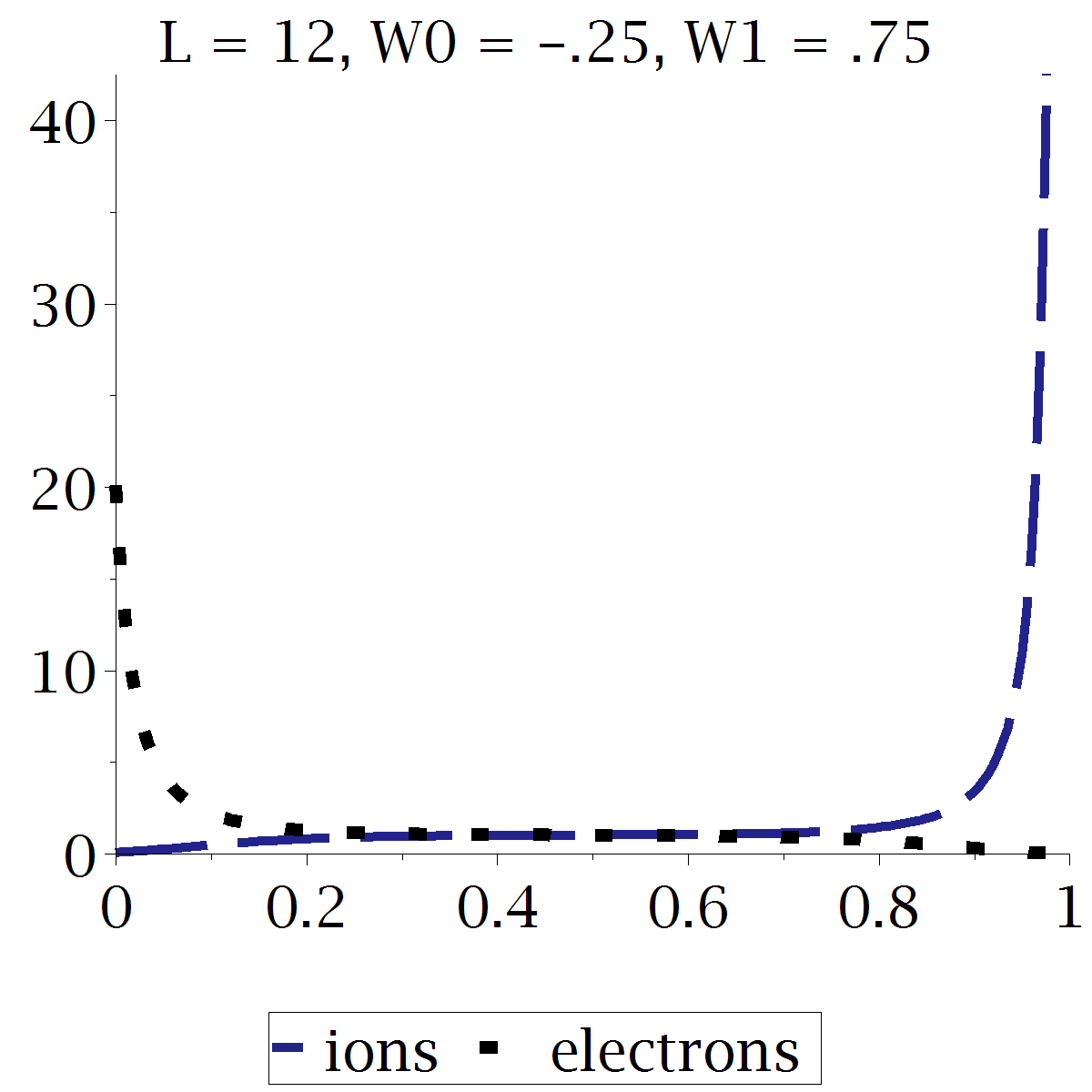}
\end{minipage}\hskip0ex 
\begin{minipage}[b]{.32\textwidth}
\includegraphics[height=.7\textwidth]{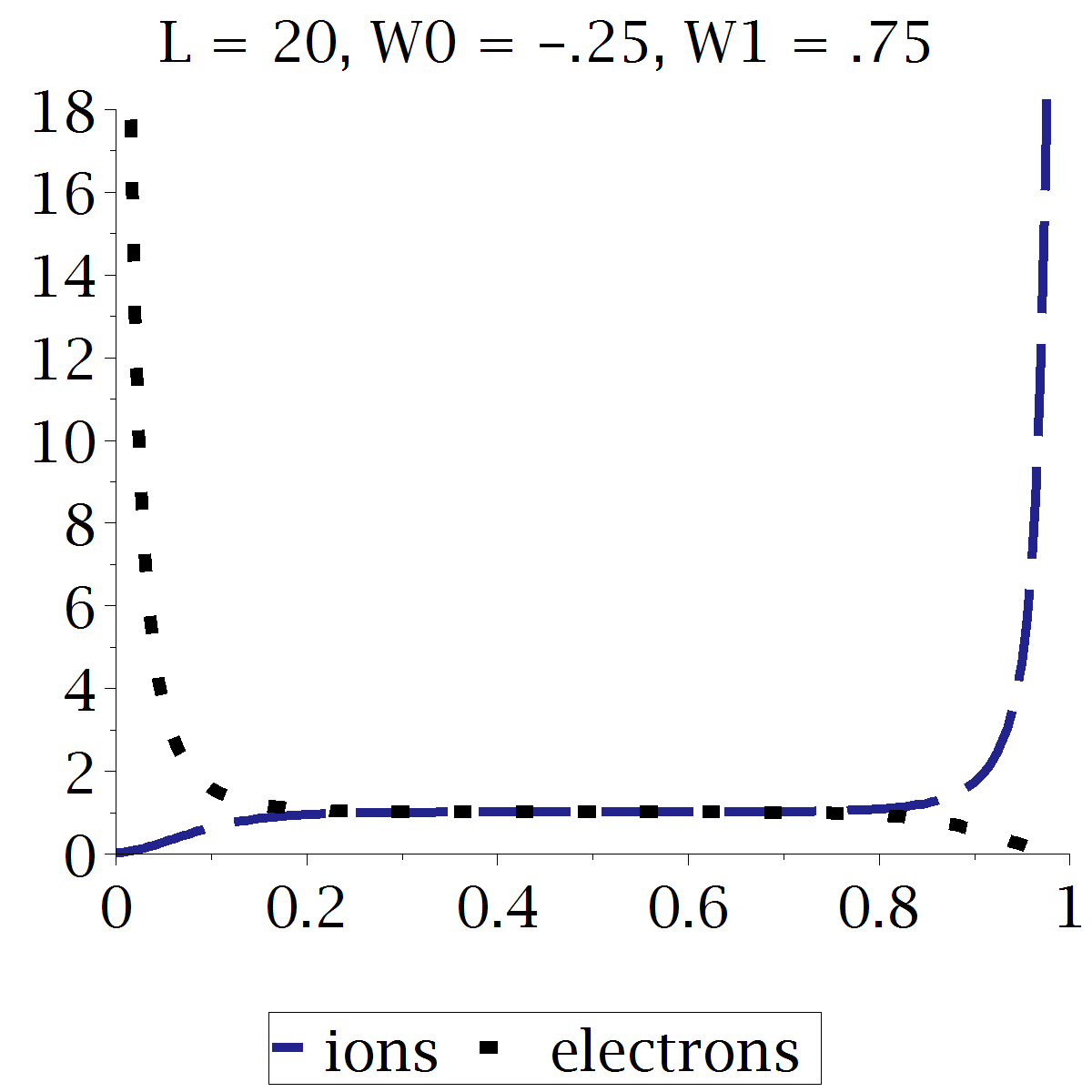}
\end{minipage}
\captionof{figure}{graphs of ion-electron density functions}\label{fig:graph-pm}
\end{minipage}

\subsubsection{approximate zero} 

The solution to $\ww(x)=0$, equivalently $\yy(x)=\qq$ (see~\eqref{eq:gentrs-v1}), represents the point where the ion/electron densities become equal. Obviously, for $\WW_1=-\WW_0$, the zero is located at $0.5$, in the middle. The explicit expression of $\ww$ allows estimating it. 

\begin{m-corollary}\label{cor:zero}
An approximate value of the zero $x_0$ of the solution $\ww_\exc$ of~\eqref{eq:gentrs-v2} is obtained by solving $\ww_1(x)=\ww_0(1-x)$; it equals $x_0\approx 0.5-\frac{2\tanh^{-1}(T)}{\LL\CC_\ave}$, with $T, \CC_\ave$ below.
\end{m-corollary}

In the boundary layer case, the values of $\CC_0$ and $\CC_1$ are about the same (almost $1$), so we replace them by $\CC_\ave:=0.5(\CC_0+\CC_1)$. 
A computation leads to a quadratic equation in $T:=\tanh\Big(\frac{\LL\CC_\ave(0.5-x)}{2}\Big)$, whose solution is: 
$$
T=\frac{\CC_\ave(\re^{\LL\WW_1/2}-\re^{\LL|\WW_0|/2})}
{\re^{\frac{\LL(\WW_1+|\WW_0|)}{2}}-\CC_\ave^2+
\sqrt{\Big(\re^{\frac{\LL(\WW_1+|\WW_0|)}{2}}-\CC_\ave^2\Big)^2
-\CC_\ave^2\big(\re^{\LL\WW_1/2}-\re^{\LL|\WW_0|/2}\big)^2}}.
$$
We observe that, as $\LL$ increases, $x_0$ drifts towards the middle. 
\begin{center}\renewcommand{\arraystretch}{1.1}
\begin{tabular}{|c|c|c|c|c|}
\hline 
$\LL$&7&12&20&30
\\\hline 
$x_0=\;$&0.44838&0.48205&0.49591&0.49921
\\\hline
$x_0\approx\;$&0.44696&0.48201&0.49591&0.49921
\\\hline 
\end{tabular}\captionof{table}{solution of $\ww(x)=0$, for $\WW_0=-0.25, \WW_1=0.75$}
\renewcommand{\arraystretch}{1}\end{center}


\section{Conclusions}\label{sct:concl}

\nit The objective of this article is twofold. 
\begin{enumerate}[leftmargin=5ex]
\item[(i)] 
We generalize the classical two-point Troesch-{\bvp} in order to take into account charge separation, which is required by physical considerations:
\begin{align*}
\ww''=\LL\sinh(\LL\ww),\;\ww(0)=\WW_0<0,\;\ww(1)=\WW_1>0\;\;(\LL>0).
\tag{$\ast$}
\end{align*}
Its solution possesses two boundary layers, for $\LL$ large enough, in contrast with the classical case. An unexpected superposition principle allows to approximately decouple $(\ast)$ into two Troesch-type equations (with modified boundary conditions). 

This leads to the main achievement: we determine an approximate, explicit solution for $(\ast)$; this is done by applying the techniques described in (ii) below. We verify the precision of the approximate solution by both numerical and analytical methods. The find that it's accurate and the precision improves for increasing $\LL$.
\item[(ii)] 
We develop a procedure for constructing explicit, approximate solutions of autonomous, second order {\bvp}s. This is important in situations when the solution possesses boundary layer(s) ---as in the case of $(*)$ above--- which make difficult using numerical algorithms. The technical approximation required by the procedure is deferred to the Appendix, to prioritize the application.
\end{enumerate}


\appendix

\section{A challenging transcendental equation}\label{sct:back}

We consider the overdetermined {\bvp}: 
\begin{align*}
\zz'=\kk\cdot(z-\alpha\CC)(\zz+\beta\CC),\;\;&\zz(0)=\zz_0,\;\zz(1)=\zz_1,
\quad(\zz_1>\zz_0>0,\;\kk>0,\;\beta, \alpha>0).
\tag{$\star\star$}\label{eq:bvp}
\end{align*}
In here, only $\alpha, \beta, \zz_0, \zz_1$ are variables. The value of $\CC\in\mbb R$ is uniquely determined by requiring the problem to be solvable. The general solution of~\eqref{eq:bvp} is: 
\begin{m-eqn}{
\frac{\zz-\alpha\CC}{\zz+\beta C}=K_1\cdot\re^{-2\LL\CC(1-x)}.
}\label{eq:bvp-sol}
\end{m-eqn}
By imposing the boundary conditions, we obtain the transcendental equation satisfied by $\CC$: 
\begin{m-eqn}{
\frac{\zz_1-\alpha\CC}{\zz_1+\beta\CC}
=
\frac{\zz_0-\alpha\CC}{\zz_0+\beta\CC}\cdot\re^{2\LL\CC},\;\;\LL:=\kk\cdot\frac{\beta+\alpha}{2}.
}\label{eq:C}
\end{m-eqn}
Thus the solution to~\eqref{eq:bvp-sol} is the following: 
\begin{m-eqn}{
K_1=\frac{\zz_1-\alpha\CC}{\zz_1+\beta\CC},\quad\zz(x)=\CC\cdot\frac{\beta K_1+\alpha\re^{2\LL\CC(1-x)}}{\re^{2\LL\CC(1-x)}-K_1}.
}\label{eq:zz}
\end{m-eqn} 
Unfortunately, it's not helpful: it's impossible to exactly solve~\eqref{eq:C} and, even numerically, it's not clear how to deal with the problem, as it contains the unknown parameter $\CC$. This compels us to approximate the solution using analytic methods. 

\subsection{Case $\alpha=\beta=1$}\label{ssct:treq-1}

We consider $\zeta>1$, $\Lda>0$, and the equation in the unknown $\cc\ges0$: 
\begin{m-eqn}{
\begin{array}{r}
\disp\frac{\zeta-\cc}{\zeta+\cc}=\frac{1-\cc}{1+\cc}\cdot e^{2\Lda \cc}
\;\Leftrightarrow\;
(\zeta-1)\cdot\frac{\cc}{\zeta-\cc^2}=\tanh(\Lda \cc)
\;\Leftrightarrow\;
\zeta= \underbrace{\cc\cdot\frac{1-\cc\cdot \tanh(\Lda \cc)}{\cc-\tanh(\Lda \cc)}}_{rhs(\cc)}.
\end{array}
}\label{eq:cc}
\end{m-eqn}
As function of $\cc$, the right-hand side has a pole at $\cc_0(\Lda)$, satisfying $\cc=\tanh(\Lda \cc)$:
\begin{itemize}[leftmargin=3ex]
\item 
for $\Lda\in(0,1]$, the solution is $\cc_0(\Lda)=0$;
\item 
for $\Lda>1$, there is a (unique) solution $\cc_0(\Lda)\in(0,1)$;
\item 
in both cases, $\cc_0(\Lda)$ can be expressed as an infinite iteration of the `tanh' function: 
$$
\cc_0(\Lda)=\tanh(\Lda\cdot\tanh(\Lda\cdot\tanh(\Lda\cdot \dots))).
$$
\end{itemize}
This formula allows approximating $\cc_0(\Lda)$. Numerical tests show that $30$ iterations of $\tanh(\Lda\cdot)$ yield at least a two digit approximation ---an upper bound--- of $\cc_0(\Lda)$, for $\Lda\ges1.1$. 
\renewcommand{\arraystretch}{1.1}\begin{m-eqn}{
\begin{array}{|l|l||l|l|}
\hline 
\multicolumn{4}{|c|}{\Lda\cc_0(\Lda)=\tanh(\Lda\cdot\tanh(\Lda\cdot\tanh(\Lda\cdot \dots)))} \\ 
\hline 
\Lda=1&\cc_0(1)=0 & \Lda=1.67&\cc_0(1.67)\approx0.9\\ 
\hline 
\Lda=1.10&\cc_0(1.10)\approx0.5&\Lda=2&\cc_0(2)\approx 0.95\\ 
\hline 
\Lda=1.50&\cc_0(1.50)\approx0.85&\Lda\ges5&\cc_0(\Lda)>0.999\approx1\\ 
\hline 
\end{array}
}\label{tab:cc}
\end{m-eqn}\renewcommand{\arraystretch}{1}

\begin{m-lemma}\label{lm:cc}
For $\zeta>1, \Lda>0$, the equation~\eqref{eq:cc} admits a unique solution $c=c(\zeta,\Lda)$: 
\begin{enumerate}
\item if $\Lda>1$,  one has $\cc\in(\cc_0(\Lda),1)$.
\item for $\Lda\in(0,1]$, $\zeta<\frac{1}{1-\Lda}$, one has $\cc\in(0,1)$.  
(For $\Lda=1$, we let $1/(1-1)=+\infty$.)
\end{enumerate}
\end{m-lemma}

\begin{proof}
In the first case, we have $\uset{\cc\to \cc_0(\Lda)^+}{\lim}rhs(\cc) =+\infty$ and $rhs(1)=1<\zeta$. In the second case, the denominator doesn't vanish, and we have $\uset{\cc\to 0^+}{\lim}rhs(\cc) =\frac{1}{1-\Lda}>\zeta$, $rhs(1)=1<\zeta$. Since the function $c\mt rhs(c)$ is decreasing, the uniqueness follows. 
\end{proof}



Now we explicitly determine upper and lower bounds for the exact solution of~\eqref{eq:cc}. Lemma~\ref{lm:cc} and Table~\ref{tab:cc} show that we should distinguish between two cases: 
\begin{itemize}[leftmargin=3ex]
\item 
$\Lda\ges1.5$, when the value of $\cc(\zeta,\Lda)$ is close to $1$;
\item 
$\Lda<1.5$, $\cc(\zeta,\Lda)$ approaches $0$, and a different approach is needed. 
\end{itemize}


\subsubsection{Case $\Lda\ges1.5$} 

Polynomial truncations fail because the $\tanh$-function approaches exponentially fast the value $1$. Since $\cc$ belongs to $(0,1)$, we write: 
$$\;\cc=\tanh(Y)=\tanh(\Lda X)>\cc_0(\Lda)=\tanh(\Lda\cdot \cc_0(\Lda)),\; X>\cc_0(\Lda).$$ 

The addition formula for the $\tanh$-function yields:
\begin{m-eqn}{
\zeta=\frac{\cc}{\tanh(Y-\Lda \tanh(Y))}=\frac{\tanh(\Lda X)}{\tanh(\Lda(X-\tanh(\Lda X)))},\quad X=\cc_0(\Lda)+\gamma,\;\gamma>0.
}\label{eq:X}
\end{m-eqn}
To determine explicit approximations, we consider the inequalities `$\les$' and `$\ges$'. 

\begin{itemize}
\item[`$\les$'] 
Since $\cc>\cc_0(\Lda)$, it's enough to have $\tanh\big(\Lda (X-\tanh(\Lda X))\big)\les\zeta^{-1}\cc_0(\Lda)$, so\\ 
$X-\tanh(\Lda X)\les \Lda^{-1}\cdot\tanh^{-1}(\zeta^{-1}\cdot \cc_0(\Lda)).$
The left-hand side is at most\\ 
$\gamma+\cc_0(\Lda)-\tanh(\Lda \cc_0(\Lda))=\gamma$. Thus the following value is convenient:
\begin{m-eqn}{
\cc_-(\zeta,\Lda)=\tanh\biggl( \Lda\cdot \cc_0(\Lda) + \tanh^{-1}(\zeta^{-1}\cdot \cc_0(\Lda)) \biggr). 
}\label{eq:ccm}
\end{m-eqn}
\item[`$\ges$'] 
As $c<1$, it's enough to have $\tanh\big(\Lda (X-\tanh(\Lda X))\big)\ges\zeta^{-1}\ges\zeta^{-1}\cc$, so 
$$
\begin{array}{l}
\gamma+\tanh(\Lda \cc_0(\Lda))-\tanh(\Lda X)=X-\tanh(\Lda X)\ges\Lda^{-1}\tanh^{-1}(\zeta^{-1}).
\end{array}
$$
A concavity argument yields $\tanh(\Lda X)-\tanh(\Lda \cc_0(\Lda))\les\frac{\Lda\gamma}{\cosh^2(\Lda \cc_0(\Lda))}$, which implies that the following value is convenient:
\begin{m-eqn}{
\cc_+(\zeta,\Lda)= \tanh\biggl(\Lda\cdot \cc_0(\Lda) + \frac{\tanh^{-1}(\zeta^{-1})}{1-\frac{\Lda}{\cosh^2(\Lda\cdot \cc_0(\Lda))}}\biggr).
}\label{eq:ccp}
\end{m-eqn}
By differentiating $\cc_0(\Lda)=\tanh(\Lda \cc_0(\Lda))$, one finds the denominator above is positive.
\end{itemize}


\subsubsection{Case $\Lda\ges5,\;\zeta\les e^{\Lda}$}

Table~\ref{tab:cc} shows $1-\cc_0(\Lda)<10^{-3}$, so we set $X=1+\gamma$ in~\eqref{eq:X}. The very same approach yields: 
\begin{m-eqn}{
\begin{array}{rl}
\text{`$\les$' is satisfied by}&\disp 
c_-(\zeta,\Lda):=\tanh\Big(\Lda\cdot\tanh(\Lda)+\tanh^{-1}\big(\zeta^{-1}\cdot\tanh(\Lda) \big) \Big);
\\[1ex] 
\text{`$\ges$' is satisfied by}&\disp 
c_+(\zeta,\Lda):=\tanh\Biggl(
\frac{\tanh^{-1}(\zeta^{-1})+\frac{\Lda(\sinh(2\Lda)-2\Lda)}{2\cosh^2(\Lda)}}
{1-\frac{\Lda}{\cosh^2(\Lda)}}  \Biggr).
\end{array}
}\label{eq:cc5}
\end{m-eqn}
The condition $\zeta\les e^\Lda$ ensures that $\gamma\ges0$. 


\subsubsection{Case $\Lda\in(0,1.5)$} 

In this situation, the estimates $c_0(\Lda)< c< 1$ are loose. Instead, we consider the power series expansion of $\frac{\tanh(x)}{x}$ and deduce: 
$$\Bigl|\frac{\tanh(x)}{x}-(1-0.2x^2)\Bigr|\les 0.15x^2,\,x\in[0,1.5].$$
Thus the solutions of the following two, explicitly solvable, quadratic equations 
\begin{m-eqn}{
\begin{array}{l}
\frac{\zeta-1}{\Lda}=(\zeta-\cc^2)\cdot\bigl(1 - 0.35\cdot\Lda^2\cc^2 \bigr),\quad 
\frac{\zeta-1}{\Lda}=(\zeta-\cc^2)\cdot\bigl(1 - 0.05\cdot\Lda^2\cc^2 \bigr),
\\[1ex]\disp
\cc_\rho(\zeta,\Lda)=
{\Biggl[ \frac{2\frac{1-\zeta(\Lda-1)}{\Lda}}{(1+\rho\cdot\zeta\Lda^2)+\sqrt{(1-\rho\cdot\zeta\Lda^2)^2+4\rho\Lda(\zeta-1)}}
\Biggr]}^{1/2},\; \rho=0.35\;and\;0.05,
\end{array}
}\label{eq:c-Lsmall}
\end{m-eqn}
bound the solution of~\eqref{eq:cc}: $\,\cc_{0.35}(\zeta,\Lda)<\cc<\cc_{0.05}(\zeta,\Lda).$ 
We insert into~\eqref{eq:X} and obtain: 
\begin{itemize}
\item[`$\les$'] $\tanh\big(\Lda (X-\tanh(\Lda X))\big)\approx\zeta^{-1}\cc_{0.35}$ yields: 
$$c_-(\zeta,\Lda):=\tanh\Big(\Lda\cdot\cc_{0.35}+\tanh^{-1}\big(\zeta^{-1}\cdot\cc_{0.35}\big)\Big).$$
\item[`$\ges$'] $\tanh\big(\Lda (X-\tanh(\Lda X))\big)\approx\zeta^{-1}\cc_{0.05}$ yields: 
$$c_+(\zeta,\Lda):=
\tanh\Biggl(
\frac{\tanh^{-1}(\zeta^{-1}\cc_{0.05})+\frac{\Lda\big(\sinh(2\Lda\cc_{0.05})-2\Lda\cc_{0.05}\big)}{2\cosh^2(\Lda\cc_{0.05})}}
{1-\frac{\Lda\cc_{0.05}}{\cosh^2(\Lda\cc_{0.05})}}  \Biggr).$$
\end{itemize}


\begin{m-remark}\label{rmk:cc} 
\begin{enumerate}[leftmargin=5ex]
\item 
Note that $\cc_+(\zeta,\Lda)\les c\les c_-(\zeta,\Lda)$. The signs are chosen in such a way that $-$ (resp. $+$) determines lower (resp. upper) bounds for the exact solution. 
\item 
The equations~\eqref{eq:ccm},~\eqref{eq:ccp} involve $\cc_0(\Lda)$. For numerical computations, we replace it by $30$ iterations of $\tanh(\Lda\cdot)$: 
$\,c_0(\Lda)\approx\underbrace{\tanh(\Lda\cdot\tanh(\dots(\tanh(\Lda)\dots)))}_{30\;\text{times}},\;\Lda\in[1.5,5].$
\item[] 
For values $\Lda\ges5$, we use the estimates~\eqref{eq:cc5}. 
\end{enumerate}
\label{rmk:treq-2}
\end{m-remark}


\subsection{Case $\alpha, \beta$ arbitrary} \label{ssct:treq-2}

For $z_1>z_0>0$, $\LL>0$, $\bta\ges\ala>0$, we consider 
\begin{m-eqn}{
\frac{z_1-\ala\CC}{z_1+\bta\CC}=\frac{z_0-\ala\CC}{z_0+\bta\CC}\cdot e^{2\LL\CC}, 
}\label{eq:zetac}
\end{m-eqn}
with $\CC\ges0$ unknown. The substitutions 
$$
\zeta_i:=z_i+\frac{\bta-\ala}{2}\CC\;(i=0,1),\quad\cc:=\frac{\bta+\ala}{2\zeta_0}\CC,\quad\Lda:=\frac{2\zeta_0\LL}{\bta+\ala}, \quad \zeta:=\frac{\zeta_1}{\zeta_0},
$$ 
reduce this equation to~\eqref{eq:cc}, and $\CC$ is given by 
\begin{m-eqn}{
\CC=\frac{2z_0\cc}{(\bta+\ala)-(\bta-\ala)\cc}.
}\label{eq:CC}
\end{m-eqn}
Complications occur, as $\zeta,\cc$ are related.  (Before, they were independent.) Fortunately, we need solving the inequalities `$\les$' and `$\ges$' in~\eqref{eq:zetac}, which leads to the same inequalities for 
$$
\zeta=\frac{z_1+\frac{\bta-\ala}{2}\CC}{z_0+\frac{\bta-\ala}{2}\CC}
=\cc\cdot\frac{1-\cc\cdot\tanh(\Lda\cc)}{\cc-\tanh(\Lda\cc)}. 
$$

\begin{itemize}
\item The left-hand side is decreasing in $\CC$ (see~\eqref{eq:CC}), and we have 
\begin{m-eqn}{
\frac{z_0c_0(\Lda)}{\beta}\les\frac{2z_0\cc_0(\Lda)}{(\bta+\ala)-(\bta-\ala)\cc_0(\Lda)}<\CC<\frac{z_0}{\ala}.
}\label{eq:CC-bounds}
\end{m-eqn}
\item The right-hand side is  increasing in $\Lda$, and we have $\frac{z_0}{\bta}L\les\Lda\les\frac{z_0}{\ala}L$.
\end{itemize}
We conclude the following:
\begin{itemize}
\item[`$\les$'] Requires 
$\zeta_-:=\frac{z_1+\frac{(\bta-\ala)z_0}{2\bta}c_0(\frac{z_0}{\beta}L)}{z_0+\frac{(\bta-\ala)z_0}{2\bta}c_0(\frac{z_0}{\beta}L)}\les 
\cc\cdot\frac{1-\cc\cdot\tanh(\frac{z_0}{\bta}L\cc)}{\cc-\tanh(\frac{z_0}{\bta}L\cc)},$ 
satisfied by $\cc=\cc_-(\zeta_-,\frac{z_0}{\bta}L)$. 
\item[`$\ges$'] Requires 
$\zeta_+:=\frac{z_1+\frac{(\bta-\ala)z_0}{2\ala}}{z_0+\frac{(\bta-\ala)z_0}{2\ala}}\ges
\cc\cdot\frac{1-\cc\cdot\tanh(\frac{z_0}{\ala}L\cc)}{\cc-\tanh(\frac{z_0}{\ala}L\cc)},$ 
satisfied by $\cc=\cc_+(\zeta_+,\frac{z_0}{\ala}L)$. 
\end{itemize}
The corresponding values $C_\pm$ are obtained by inserting $c_\pm$ into~\eqref{eq:CC}.



\begin{thebibliography}{ooo}

\bibitem{br-so} \bibauth{Braga}{F.} \bibauth{Soares}{D.} \bibtitl{Linear pinch equilibrium of non-neutral plasma revisited:phenomenological consequences of a numerical accuracy problem} \bibjnyp{Plasma Phys. Tech.}{6}{2019}{217-222}

\bibitem{fi-tr} 
\bibauth{Firnett}{P.} \bibauth{Troesch}{B.} \bibtitl{Shooting-splitting method for sensitive two-point boundary value problems} 
\bibinbook{Proc. Conference on the Numerical Solution of Ordinary Differential Equations. D.~Bettis (ed.)}{Lect. Notes Math. 362}{Springer, Berlin, Heidelberg, 1974}{pp. 408-433} 

\bibitem{csfx}
\bibauth{Halic}{M.} \bibauth{Tajarod}{R.} \bibauth{Temimi}{H.} \bibtitl{An Accurate Analytical Solution to Troesch’s Problem Through Lower
and Upper Envelope Techniques} \bibjnyp{Chaos, Solitons \& Fractals: X}{11}{2023}{100096}

\bibitem{so-br} \bibauth{Soares}{D.} \bibauth{Belich}{H.} \bibauth{Spalenza}{W.} \bibauth{Braga}{F.} \bibtitl{Nonneutral Weibel model plasma in the non-minimal CPT-odd coupling} \bibjnyp{Eur. Phys. J. C}{84}{2024}{232}

\bibitem{tsch} \bibauth{Troesch}{B.} \bibtitl{A Simple Approach to a Sensitive Two-Point Boundary Value Problem} \bibjnyp{J. Comp. Phys.}{21}{1976}{279-290} 

\bibitem{weib1} \bibauth{Weibel}{E.} \bibtitl{The Plasma in a Magnetic Field} \bibinbook{Symposium Magnetohydrodynamics. R. Landshoff (ed.)}{Stanford Univ. Press}{1958}{pp. 60-76}

\bibitem{weib2} \bibauth{Weibel}{E.} \bibtitl{On the Confinement of a Plasma by Magnetostatic Fields} \bibjnyp{Phys. Fluids}{2}{1959}{52-56}

\bibitem{ya-kr} \bibauth{Yamakawa}{T.} \bibauth{Kreinovich}{V.} \bibtitl{Why Fundamental Physical Equations Are of Second Order} \bibjnyp{Internat. J. Theor. Phys.}{38}{1999}{1763-1769}
\end{thebibliography}
\end{document}